\documentclass[11pt]{article}
\usepackage{a4wide}

\usepackage{graphicx,amssymb,amsmath}
\usepackage{epsfig}
\usepackage{epsf}
\usepackage{latexsym}
\usepackage{amssymb}
\usepackage{amsthm}
\usepackage{float}

\usepackage{color}
\usepackage{hyperref} 
\definecolor{modra3}{rgb}{.1,.0,.4}
\hypersetup{colorlinks=true, linkcolor=modra3, urlcolor=modra3, citecolor=modra3, pdfpagemode=UseNone, bookmarks=true, pdfstartview=}

\newtheorem{theorem}{Theorem}
\newtheorem{lemma}[theorem]{Lemma}
\newtheorem{problem}{Problem}

\newtheorem{observation}[theorem]{Observation}

\newcommand{\old}[1]{{}}

\def\P{\mathcal P}

\def\cro{{\mbox {\sc cr}}}

\begin{document}

\title{Saturated simple and $k$-simple topological graphs}

\author{%
Jan Kyn\v cl\thanks{Department of Applied Mathematics and Institute for Theoretical Computer Science,
Charles University, Faculty of Mathematics and Physics, 
Malostransk\'e n\'am.~25, 118~00~ Praha 1, Czech Republic; and 
Alfr\'ed R\'enyi Institute of Mathematics, Budapest, Hungary.
Email: \texttt{kyncl@kam.mff.cuni.cz}. 
Supported 
by the GraDR EUROGIGA GA\v{C}R project No. GIG/11/E023, by the grant
SVV-2013-267313 (Discrete Models and Algorithms) and by ERC Advanced Research
Grant no 267165 (DISCONV).}
\and 
J\'anos Pach\thanks{Ecole Polytechnique F\'ed\'erale de Lausanne and
Alfr\'ed R\'enyi Institute of Mathematics, Budapest. Email: \texttt{pach@cims.nyu.edu}.
Research partially supported by Swiss National Science Foundation Grants
200021-137574 and 200020-144531, by Hungarian Science Foundation Grant OTKA NN
102029 under the EuroGIGA programs ComPoSe and GraDR, and by NSF grant CCF-08-30272.}
\and
Rado\v{s} Radoi\v{c}i\'{c}\thanks{Department of Mathematics, Baruch College, City
University of New York, NY, USA.
Email: \texttt{rados.radoicic@baruch.cuny.edu}.
Part of the research by this author was done at the Alfr\'ed R\'enyi Institute of
Mathematics, partially supported by Hungarian Science Foundation Grant OTKA T 046246.}
\and
G\'eza  T\'oth\thanks{Alfr\'ed R\'enyi Institute of Mathematics,
Budapest, Hungary. Email: \texttt{geza@renyi.hu}.
Partially supported by Hungarian Science Foundation Grants OTKA  K 83767 and NN 102029.}}

\maketitle

\begin{abstract} A {\em simple topological graph} $G$ is a graph drawn in the plane so that any
pair of edges have at most one point in common, which is either
an endpoint or a proper crossing. $G$ is called {\em saturated\/} if no further edge can be added
without violating this condition.
We construct saturated simple topological graphs with $n$ vertices and $O(n)$ edges.
For every $k>1$, we give similar constructions for  
{\em $k$-simple topological graphs}, that is,  for graphs drawn in the plane so that any
two edges have at most $k$ points in common.
We show that in any $k$-simple topological graph, any two independent
vertices can be connected by a curve that crosses each of the original edges
at most  $2k$ times. Another construction shows  that the bound $2k$ cannot be improved.
Several other related
problems are also considered. 
\end{abstract}


\section{Introduction}

Saturation problems in graph theory have been studied at length, ever since
the paper of Erd\H{o}s, Hajnal, and Moon~\cite{EHM64}. Given a
graph $H$, a graph $G$ is $H$-saturated if $G$ does not contain $H$ as a subgraph,
but the addition of any edge joining two non-adjacent vertices of $G$
creates a copy of $H$. The saturation number of $H$, $\mbox{sat}(n,H)$, is the
minimum number of edges in an $H$-saturated graph on $n$ vertices.
The saturation number for complete graphs was determined in~\cite{EHM64}.
A systematic study by K\'{a}szonyi and Tuza~\cite{KT86} found the best known general upper
bound for $\mbox{sat}(n,H)$ in terms of the independence number of $H$.
The saturation number is now known, often precisely, for many graphs;
for these results and related problems in graph theory we refer the reader to
the thorough survey of J. Faudree, R. Faudree, and Schmitt~\cite{FFS11}.
It is worth noting that $\mbox{sat}(n,H) = O(n)$,
quite unlike the Tur\'{a}n function $\mbox{ex}(n,H)$, which is often superlinear.

In this paper, we study a saturation problem for {\em drawings\/} of graphs.
In a drawing of a simple undirected graph $G$ in the plane,
every vertex is represented by a point, and every edge is represented by a
curve between the points that correspond to its endpoints.
If it does not lead to confusion, these points and curves
are also called {\em vertices\/} and {\em edges}.
We assume that in a drawing no edge passes through a vertex and no two edges
are tangent to each other.
A graph, together with its drawing, is called a {\em simple topological graph\/}
if any two edges have at most one point in common, which is either their common
endpoint or a proper crossing. In general, for any positive integer $k$, it is called 
a {\em $k$-simple topological graph\/} if any two edges have at most $k$ points
in common.
We also assume that in a $k$-simple topological graph no edge crosses itself.
Obviously, a  $1$-simple topological graph is a simple topological graph.


Our motivation partly comes from the following problem: At least how many
pairwise disjoint edges can one find in every simple topological 
graph with $n$ vertices and $m$ edges~\cite{PST03}? (Note that the 
simplicity condition is essential here, as there are complete 
topological graphs on $n$ vertices and no {\em two\/} disjoint edges,
in which every pair of edges intersect at most twice~\cite{PT10}.)
For {\em complete\/} simple topological graphs, i.e., when $m={n\choose 2}$, 
Pach and T\'{o}th conjectured (\cite{BMP05}, page 398) that one can 
always find $\Omega(n^{\delta})$ disjoint edges for a suitable 
constant $\delta > 0$. 
This was shown by Suk \cite{S13} with $\delta=1/3$; see \cite{FR13} 
for an alternative proof. Recently, Ruiz-Vargas \cite{R13}
has improved this 
bound to $\Omega\left(\sqrt{{n}/{\log{n}}}\right)$. 
Unfortunately, all known proofs break down for {\em non-complete\/} simple 
topological graphs. For {\em dense\/} graphs, i.e., when 
$m\ge \varepsilon n^2$ for some $\varepsilon>0$,
Fox and Sudakov \cite{FS09} established the existence of  
$\Omega(\log^{1+\gamma}{n})$ pairwise disjoint edges, with 
$\gamma \approx 1/50$. However, if $m\ll n^2$, the best known
lower bound, due to Pach and T\'oth \cite{PT10}, is only
$\Omega\left((\log{m}-\log{n})/\log{\log{n}}\right)$. 

We know a great deal about the structure of complete simple
topological graphs, but in the non-complete case our knowledge
is rather lacunary. We may try to extend a simple topological
graph to a complete one by adding extra edges and then explore
the structural information we have for complete graphs.
The results in the present note suggest that this 
approach is not likely to succeed:  there exist very sparse 
simple topological graphs to which no edge can be added without 
violating the simplicity condition.


A $k$-simple, non-complete topological graph is {\em saturated\/} if no further edge can be added
so that the resulting drawing is still a $k$-simple topological graph.
In other words, if we connect any two non-adjacent vertices by a curve, it
will have at least $k+1$ common points with one of the existing edges.

\begin{figure}[ht]
\begin{center}
\includegraphics{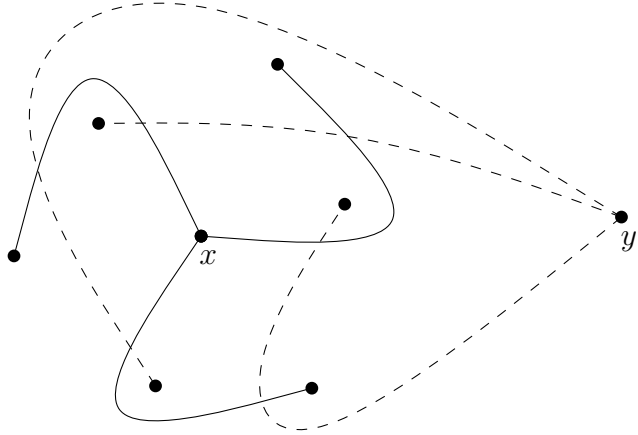}
\caption{A topological graph $G_1$: the edge $\{x, y\}$ cannot be added.\label{kampok}}
\end{center}
\end{figure}

Consider the simple topological graph $G_1$ with eight vertices, depicted
in Figure~\ref{kampok}. It is easy to verify that the vertices $x$ and $y$
cannot be joined by a new edge so that the resulting topological graph remains simple.
Indeed, every edge of $G_1$ is incident either to $x$ or to $y$,
and any curve joining $x$ and $y$ must cross at least one edge.
On the other hand, $G_1$ can be extended to a (saturated) simple topological
graph in which every pair of vertices except $x$ and $y$ are connected by an edge.

\begin{figure}[ht]
\begin{center}
\includegraphics{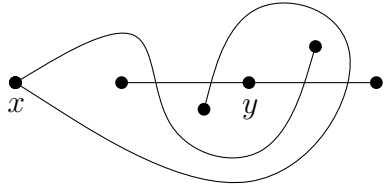}
\caption{A topological graph $G_2$: the edge $\{x, y\}$ cannot be added.\label{kynclpelda}}
\end{center}
\end{figure}

Another example was found independently by Kyn\v{c}l~\cite[Fig.~9]{K13}: The
simple topological graph $G_2$ depicted in Figure~\ref{kynclpelda}
has only six vertices, from which $x$ and $y$ cannot be joined by an edge
without intersecting one of the original edges at least twice. Again, $G_2$ can
be extended to a simple topological graph in which every pair of vertices except
$x$ and $y$ are connected by an edge. 

In view of the fact that the graphs shown in
Figures~\ref{kampok} and~\ref{kynclpelda} can be extended to nearly complete simple
topological graphs, it is a natural question to ask whether every saturated simple topological graph with $n$ vertices must have $\Omega(n^2)$ edges. 
It is not obvious at all, whether there exist saturated non-complete
$k$-simple topological
graphs for some $k>1$. 
Our next theorem shows that there are such graphs, for every $k$, moreover, they may have
only a linear number of edges.



\begin{theorem}\label{linearbound}
For any positive integers $k$ and $n\ge 4$, 
let $s_k(n)$ be the minimum number of edges that a saturated $k$-simple topological
graph on $n$ vertices can have. Then 
\begin{enumerate}
\item[(i)] we have 
$$1.5n\le s_1(n)\le 17.5n,$$
\item[(ii)] for $k>1$ we have 
$$n\le s_k(n)\le 16n.$$
\end{enumerate}
\end{theorem}

For our best upper bounds see Table~\ref{table1}.

\begin{table}
\begin{center}
{\footnotesize
\begin{tabular}{l|c|c|c|c|c|c|c|c|c|c|c}
$k$ & 1 & 2 & 3 & 4 & 5 & 6 & 7 & 8 & 9 & 10 & $\ge 11$ \\
\hline
{\mbox{upper bound}} & $17.5n$ & $16n$ & $14.5n$ & $13.5n$ & $13n$ & $9.5n$ & $10n$ & $9.5n$ & $7n$ & $9.5n$ & $7n$ 
\end{tabular}
\caption{Upper bounds on the minimum number of edges in saturated $k$-simple
  topological graphs.}
}
\label{table1}
\end{center}
\end{table}

For any positive integers $k$ and $l$, $k < l$, a topological graph $G$ {\em
together with\/}
a pair of non-adjacent vertices $\{u, v\}$ is called a 
{\em $(k,l)$-construction\/} if
$G$ is $k$-simple and any curve joining $u$ and $v$ has at least 
$l$ points in common
with at least one edge of $G$. Using this terminology, 
every saturated non-complete $k$-simple
topological graph together with any pair of non-adjacent vertices is a 
$(k,k+1)$-construction.

\begin{theorem}\label{klconstruction}
For every $k>0$,
\begin{enumerate}
\item[(i)] There exists a $(k, 2k)$-construction,
\item[(ii)] There is no $(k, l)$-construction with $l > 2k$.
\end{enumerate}
\end{theorem}

For any positive integers $k$ and $l$, $k < l$,
a {\em non-complete\/} topological graph $G$ is called {\em $(k, l)$-saturated\/} if
$G$ is $k$-simple and any curve joining {\em any pair\/} of non-adjacent
vertices has at least $l$ points in common with at least one edge of $G$.
Obviously, every saturated $k$-simple topological graph is $(k,
k+1)$-saturated.
Clearly, every $(k, l)$-saturated topological graph, together with any pair of its
non-adjacent vertices, is a $(k, l)$-construction. 
However, for $l>k+1$, the existence of a $(k, l)$-construction
does not necessarily imply the existence of a $(k, l)$-saturated topological graph.
The best we could prove is the following.

\begin{theorem}\label{3k/2saturated}
For any $k>0$, there exists a $(k, \lceil 3k/2\rceil)$-saturated topological graph.
\end{theorem}

In the proof of Theorem~\ref{klconstruction} we obtain a set of six curves,
any two of which cross at most once, and two points, such that any curve
connecting them has to cross one of the six curves at least twice (see Figure~\ref{nembovitheto1}).

On the contrary, it follows from Levi's enlargement
lemma~\cite{L26} that 
if the curves have to be {\em bi-infinite}, that is, two-way unbounded, 
then there
is no such construction. A {\em pseudoline arrangement\/} is a set of 
bi-infinite curves such that any two of them cross exactly 
once. 
By Levi's lemma, for any two points not on the same line, the
arrangement can be extended by a pseudoline through these two points. 
A {\em $k$-pseudoline arrangement\/} is a set of 
bi-infinite curves such that any two of them cross at most $k$ times. 
A {\em $k$-pseudocircle arrangement\/} is a set of 
closed curves such that any two of them cross at most $k$ times.
Elements of pseudoline arrangements and $k$-pseudoline arrangements are called {\em pseudolines}.
Note that for even $k$, 
$k$-pseudoline arrangements can be considered a special case of
$k$-pseudocircle arrangements. 

Snoeyink and Hershberger
\cite{SH91} 
generalized Levi's lemma to $2$-pseudoline arrangements
and  $2$-pseudocircle arrangements
as follows. 
They proved that for every $2$-pseudocircle arrangement and three points, not
all
on the same pseudocircle, the 
arrangement can be extended by a closed curve through these three points so
that
it remains a $2$-pseudocircle arrangement.
They also showed that for $k\ge 3$, an analogous statement with 
$k$-pseudoline arrangements and $k+1$ given
points is false.

A $k$-pseudoline arrangement 
is {\em $(p,l)$-forcing} if there is a set $A$ of $p$ points 
such that every bi-infinite curve through the points of $A$ crosses one 
of the pseudolines at least $l$ times. 
Snoeyink and Hershberger~\cite{SH91} found $(k+1,k+1)$-forcing $k$-pseudoline
arrangements for $k\ge 3$.
We generalize their result as follows.

\begin{theorem}\label{pseudoline}
\begin{itemize}
\item[(i)] For every $k\ge 1$, there is a $(3,5\lceil(k-7)/4\rceil)$-forcing $k$-pseudoline
arrangement.
\item[(ii)] For every $k\ge 1$, there is a $(k,\Omega(k\log k))$-forcing $k$-pseudoline arrangement.
\end{itemize}
\end{theorem}

In Section 2 we define tools necessary for our constructions. In Section 2.1 we define {\em spirals\/} and use them in Lemma~\ref{lemma_weakklconstruction} to prove 
the existence of a $\left( k, \lceil 7(k-1)/6\rceil \right)$-construction (for $k\ge 8$).
Although Lemma~\ref{lemma_weakklconstruction} is a very weak version of Theorem~\ref{klconstruction},
its proof is a relatively simple construction, which serves as the basis of all our further constructions.
In Section 2.2 we define another tool, {\em forcing blocks}, and as an illustration, 
we prove Lemma~\ref{2k-1construction}, which is an improvement of 
Lemma~\ref{lemma_weakklconstruction}, yet still weaker than Theorem~\ref{klconstruction}.
Finally, in Section 2.3, we define the remaining necessary tools, {\em grid blocks} and 
{\em double-$k$-forcing blocks}, and use them to prove Theorem~\ref{klconstruction} (i). 

In Section 3 we prove Theorem~\ref{klconstruction} (ii). Our proof is self-contained
and independent of the tools developed in Section 2. 
In Section 4 we prove the upper bounds in Theorem~\ref{linearbound}, by giving 
five different constructions; the first one is for $k=1$, and it is essentially different from the other four, 
which use spirals, grid blocks, and forcing blocks described previously in Section 2. 
In Section 5 we prove the lower bounds in Theorem~\ref{linearbound}. Our proof is self-contained 
and independent of the remaining sections. 
In Section 6, we prove Theorem~\ref{3k/2saturated}. Our construction uses grid blocks and double-forcing 
blocks described in Section 2. 
In Section 7, we prove Theorem~\ref{pseudoline}. Our constructions use spirals from Section 2.
We finish the paper with some remarks and open problems.


\section{Building blocks for $(k, l)$-constructions}

With the exception of the proof of Theorem~\ref{linearbound} (i), 
we construct drawings on a vertical 
cylinder, which can be transformed into a planar drawing.
The cylinder will be represented by an axis-parallel rectangle whose left and
right vertical sides are identified. Curves on the cylinder are also represented 
in the  axis-parallel rectangle, where they can ``jump'' between the left and the
right sides.
Edges will be 
drawn as $y$-monotone curves.

Drawings will be constructed from {\em blocks}. Each 
block is a horizontal ``slice'' of the cylinder, represented again by an
axis-parallel rectangle, say, $R$,  whose left and
right vertical sides are identified. 
A {\em cable} in a block is a group of intervals of edges that go very close
to each other but do not cross in the block. A cable is represented by a single 
curve which goes very 
close to each edge in the cable. 
A curve or a cable in a block $B$ whose endpoints are on the top and
bottom boundary of $R$ is called a {\em transversal\/} of $B$.
For any curves or cables $a$ and $b$, let $\cro(a,b)$ denote the
number of crossings between $a$ and $b$.

\subsection{\texorpdfstring{Spirals and a $\left(k, \lceil 7(k-1)/6\rceil \right)$-construction}{Spirals and a $(k, 7(k-1)/6)$-construction}}

\begin{figure}
\begin{center}
\epsfig{file={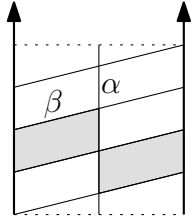}}
\end{center}
\caption{A $4$-spiral formed by $\alpha$ and $\beta$. 
One of the four regions is shaded. The arrows represent the two edges 
of the rectangle that are glued together to form a cylinder.}
\label{figure_spiral}
\end{figure}

Let $B$ be a block on the cylinder, represented by the unit square $R$ 
with vertices $(0,0)$, $(1,0)$, $(1,1)$, $(0,1)$, the two vertical sides are
identified.
Let $\alpha$ be a straight line segment from $(0.5,0)$ to $(0.5,1)$, and let
$\beta$ be 
represented as the union of $m$ straight line segments, 
$b_1, b_2, \ldots , b_m$, where $b_i$ is the segment from $(0, (i-1)/m)$ to
$(1,i/m)$. See Figure~\ref{figure_spiral}.
Cables $a$ and $b$ in $B$ form an {\em $m$-spiral\/} if there is a homeomorphism of $B$ 
that takes $a$ to $\alpha$ and $b$ to $\beta$ and maps the lower boundary of $B$ to itself. 
Clearly, such cables $a$ and $b$ are transversals of $B$ and intersect exactly $m$
times.

\begin{observation}\label{obs_spiral}
Suppose that $a$ and $b$ form an $m$-spiral in
a block $B$. 
Then every transversal of $B$ crosses $a$ and $b$
together at least $m-1$ times. 
\end{observation}

\begin{proof}
Let $\kappa$ be a transversal of $B$. 
Extend $B$ to a two-way infinite cylinder, $B'$. 
Cables $a$ and $b$ together divide the cylinder $B'$ into $m$ regions, say, $B_1, B_2,
\ldots, B_{m}$, from bottom to top. 
One endpoint of $\kappa$ is in $B_1$, the other one is in $B_{m}$.
It is easy to see that $B_i$ and $B_j$ have a common boundary
if and only if $|i-j|=1$. Therefore, to go from $B_1$ to $B_m$,
$\kappa$ has to cross at least $m-1$ boundaries. 
\end{proof}

Using spirals, we are able to prove the following weak version of Theorem~\ref{klconstruction} (i). 

\begin{lemma}\label{lemma_weakklconstruction}
For $k\ge 8$ there exists a $(k, l)$-construction with $l=\lceil
7(k-1)/6\rceil >k$.  
\end{lemma}

\begin{proof}
The construction consists of $7$ consecutive blocks, 
$X,\allowbreak A,\allowbreak B,\allowbreak C,\allowbreak D,\allowbreak
E,\allowbreak Y$, 
in this order (say, from bottom to top). 
First we define six independent edges, 
$\alpha_1,\allowbreak \alpha_2,\allowbreak \beta_1,\allowbreak 
\beta_2,\allowbreak \gamma_1,\allowbreak \gamma_2$, and two isolated
vertices, $x$ and $y$.

Put $x$ in $X$ and $y$ in $Y$.
The edges $\alpha_1$ and $\alpha_2$ are in the blocks $A$ and $B$, both have
one endpoint on the boundary of $X$ and $A$ and one on the boundary of $B$ and
$C$. 
Edges $\beta_1$ and $\beta_2$ are in $B$, $C$ and $D$, 
both have one endpoint on the boundary of $A$ and $B$ and one on the boundary of $D$ and
$E$. 
Edges $\gamma_1$ and $\gamma_2$ are in $D$ and $E$,
both have one endpoint on the boundary of $C$ and $D$ and one on the boundary of $E$ and
$Y$.
The edges $\alpha_1$ and $\alpha_2$ form a $k$-spiral in $A$ and a cable in $B$. The edges
$\beta_1$ and $\beta_2$ form another cable in $B$, and these two cables form a
$k$-spiral. Further, $\beta_1$ and $\beta_2$ form a $k$-spiral in $C$ and a cable in $D$. 
The edges $\gamma_1$ and $\gamma_2$ form another cable in $D$, and these two cables form a
$k$-spiral. 
Finally, $\gamma_1$ and $\gamma_2$ form a $k$-spiral in $E$.

We show that every curve $\kappa$ from $x$ to $y$ crosses one of the curves 
$\alpha_1,\allowbreak \alpha_2,\allowbreak 
\beta_1,\allowbreak \beta_2,\allowbreak \gamma_1,\allowbreak \gamma_2$ at least 
$7(k-1)/6$ times. Let $\kappa$ be a fixed curve from $x$ to $y$.
For every $\chi\in\{\alpha_1,
\alpha_2, 
\beta_1, \beta_2, \gamma_1, \gamma_2\}$ and $Z\in\{A, B, C, D, E\}$, 
let $Z(\chi)$ denote the number of intersections of $\chi$ with $\kappa$ in $Z$.
(That is, $A(\alpha_1)$ is the number of intersections between $\alpha_1$ and
$\kappa$ in $A$.)
By Observation~\ref{obs_spiral}, we have
\begin{align*}
A(\alpha_1)+A(\alpha_2)&\ge k-1,\\
B(\alpha_i)+B(\beta_j)&\ge k-1, \ \ 1\le i, j\le 2,\\
C(\beta_1)+C(\beta_2)&\ge k-1,\\
D(\beta_i)+D(\gamma_j)&\ge k-1, \ \ 1\le i, j\le 2,\\
E(\gamma_1)+E(\gamma_2)&\ge k-1.\\
\end{align*}
Therefore, we can assume without loss of generality
that 
$A(\alpha_1)\ge (k-1)/2$, 
$C(\beta_1)\ge (k-1)/2$ and 
$E(\gamma_1)\ge (k-1)/2$.

It follows that 
$$\cro(\alpha_1, \kappa)+
\cro(\beta_1, \kappa)+
\cro(\gamma_1, \kappa)$$
$$=A(\alpha_1)+B(\alpha_1)+
B(\beta_1)+C(\beta_1)+D(\beta_1)+
D(\gamma_1)+E(\gamma_1)\ge 7(k-1)/2.$$
Consequently, at least one of 
$\cro(\alpha_1, \kappa)$,
$\cro(\beta_1, \kappa)$,
$\cro(\gamma_1, \kappa)$ is at least 
$7(k-1)/6>k$, since $k\ge 8$.
\end{proof}

\begin{figure}
\begin{center}
\epsfig{file={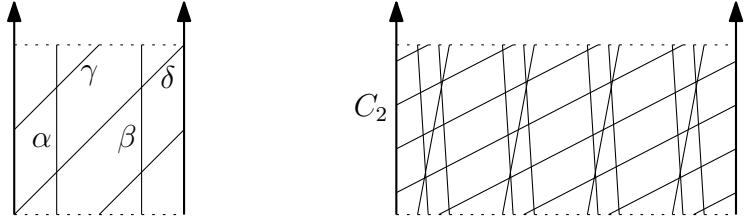}}
\end{center}
\caption{Left: a crossing-forcing configuration formed by
curves $\alpha,\beta,\gamma,\delta$. Right: cables in the subblock $C_2$.}
\label{figure_XX}
\end{figure}

\subsection{Forcing blocks and a $(k, 2k-1)$-construction}

We define another type of block, which is built
from several subblocks.
Let $B$ be a block on the cylinder represented by the unit square $R$ with vertices $(0,0)$, $(1,0)$, $(0,1)$, $(1,1)$, where the two vertical sides are identified. See Figure~\ref{figure_XX}, left. 
Let $\alpha$ be a straight line segment from $(0.25,0)$ to $(0.25,1)$, 
$\beta$ a straight line segment from $(0.75,0)$ to $(0.75,1)$,
$\delta$ a straight line segment from $(0,0)$ to $(1,1)$,
and let $\gamma$ be represented as the union of the segment
from $(0.5,0)$ to $(1,0.5)$ and the segment 
from $(0,0.5)$ to $(0.5, 1)$.
Cables $a$, $b$, $c$ and $d$  
in a block $B$ form a {\em crossing-forcing configuration\/}
if there is a homeomorphism of $B$ 
that takes $a$ to $\alpha$, 
$b$ to $\beta$, $c$ to $\gamma$, and $d$ to $\delta$, and maps the lower boundary of $B$ to itself. 
Clearly, in this case any two cables intersect at most once.

\begin{observation}\label{obs_XX}
Suppose that 
cables $a$, $b$, $c$ and $d$ form a crossing-forcing configuration
in a block $B$. Then every transversal of $B$ crosses 
at least one of $a$, $b$, $c$, or $d$. 
\qed
\end{observation}

Fix $k>0$. 
Now we 
define the {\em $k$-forcing block} $B_k$ of $4^k$ edges, $a_1,
a_2, \ldots , a_m$, $m=4^k$. 
We build $B_k$ 
from $k$ subblocks $C_1, C_2, \ldots, C_k$, 
arranged from top to bottom in this order. 
In $C_1$, divide our edges into four equal subsets, each form a cable, and
these four cables form a crossing-forcing configuration in $C_1$. 
In general, suppose that for some $i$, $1\le i<k$, $C_i$ contains $4^i$ cables 
$c_1, c_2, \ldots , c_{4^i}$ 
and 
each of them contains $4^{k-i}$ edges.
For each cable $c_j$ of $C_i$, divide the corresponding set of edges 
into four equal subsets, each of them form a cable in $C_{i+1}$, and let these
four cables form a crossing-forcing configuration in  $C_{i+1}$. 
It is possible to draw the
cables so that any two of them intersect at most once in $C_{i+1}$ and so that
for every two edges $e \in c_j$ and $f \in c_{j'}$, $j<j'$, the edges
$e$ and $f$ intersect the top and the bottom boundary of $C_{i+1}$ in the same order.
See Figure~\ref{figure_XX}, right. 
Clearly,
$C_{i+1}$ contains $4^{i+1}$ cables 
and 
each of them contains $4^{k-i-1}$ edges. 

The resulting block $B_k=\cup_{i=1}^k C_i$ 
is called a  {\em $k$-forcing block\/}
of edges $a_1, a_2, \ldots ,\allowbreak a_m$, where $m=4^k$.
The next lemma explains the name.

\begin{lemma}\label{kforcinglemma}
Suppose that $B_k=\cup_{i=1}^kC_i$ 
is a $k$-forcing block of edges 
$a_1,\allowbreak a_2,\allowbreak \ldots,\allowbreak a_m,\allowbreak m=4^k$.
Then every transversal of $B_k$ 
intersects at least one of 
$a_1,\allowbreak a_2,\allowbreak \ldots,\allowbreak a_m$ at least $k$
times.
\end{lemma}

\begin{proof}
We prove the statement by induction on $k$. 
For $k=1$ the statement is equivalent to Observation~\ref{obs_XX}.
Suppose that the statement has been proved for $k-1$, and let 
$\kappa$ be a transversal of $B_k$. 
By Observation~\ref{obs_XX}, $\kappa$ crosses at least one of the cables in $C_1$.
Consider now only the $4^{k-1}$ edges that belong to that cable.
These edges form a {\em $(k-1)$-forcing block\/} in 
$B'_{k-1}=\cup_{i=2}^kC_i$. Therefore, by the induction hypothesis, $\kappa$
crosses 
one of the edges at least $k-1$ times in $B'_{k-1}$. It also crosses this
edge in $C_1$, so we are done.
\end{proof}

Now we prove a statement which is slightly weaker than
Theorem~\ref{klconstruction} (i), but much stronger than
Lemma~\ref{lemma_weakklconstruction}.

\begin{lemma}\label{2k-1construction}
For every $k>0$,
there exists a $(k, 2k-1)$-construction.
\end{lemma}

\begin{proof}
The construction consists of $2k+1$ consecutive blocks, $X,\allowbreak F_1,\allowbreak
S_1,\allowbreak F_2,\allowbreak S_2,\allowbreak 
\ldots,\allowbreak S_{k-1},\allowbreak F_k,\allowbreak Y$,
in this order (from bottom to top). Let $m=4^k$.  
We define $km$ 
independent edges
$\alpha_i^j$, $1\le i\le k$, $1\le j\le m$, and two isolated vertices $x$ and
$y$ as follows. 
Put $x$ in $X$ and $y$ in $Y$.
\begin{itemize}
\item The edges $\alpha_1^j$, $1\le j\le m$, are in the blocks $F_1$ and $S_1$, 
\item for every $i$, $1<i<k$, the edges $\alpha_i^j$, $1\le j\le m$, are in the blocks
  $S_{i-1}$, $F_i$ and $S_i$,
\item the edges $\alpha_k^j$, $1\le j\le m$, are in the blocks $S_{k-1}$ and $F_k$.
\end{itemize}
 
For every $i$, $1\le i\le k$, the edges $\alpha_i^j$, $1\le j\le m$, form a
$k$-forcing block in $F_i$.
For every $i$, $1\le i\le k-1$, the edges  $\alpha_i^j$, $1\le j\le m$, form one
cable in $S_i$, the edges $\alpha_{i+1}^j$, $1\le j\le m$, form another cable in
$S_i$, and these two cables form a $k$-spiral in $S_i$. 

Let $\kappa$ be a fixed curve from $x$ to $y$.
We show that $\kappa$ crosses one of the curves
$\alpha_i^j$ at least $2k-1$ times. 
For every $i$, $1\le i\le k$,  by Lemma~\ref{kforcinglemma}, there is a $j$,
$1\le j\le m$, such that 
$\alpha_i^j$ and $\kappa$ cross at least $k$ times in $F_i$. Denote this 
$\alpha_i^j$ by $\alpha_i$.

For every $Z\in\{F_1, F_2, \ldots , F_k, S_1, S_2, \ldots , S_{k-1}\}$ and 
$1\le i\le k$, 
let $Z(\alpha_i)$ denote the number of intersections of 
$\alpha_i$ with $\kappa$ in the block $Z$. 
By the choice of $\alpha_i$, for every $i$, $1\le i\le k$, we have
$$F_i(\alpha_i)\ge k.$$
By Observation~\ref{obs_spiral}, for every $i$, $1\le i\le k-1$, we have
$$S_i(\alpha_i)+S_i(\alpha_{i+1})\ge k-1.$$
Summing up, 
\begin{align*}
\sum_{i=1}^k\cro(\kappa, \alpha_i)
&=\sum_{i=1}^kF_i(\alpha_i)+\sum_{i=1}^{k-1}\left(S_i(\alpha_i)+S_i(\alpha_{i+1})\right)\\
&\ge k^2+(k-1)^2=(2k-2)k+1.
\end{align*}
Therefore, for some $i$, $\cro(\kappa, \alpha_i)\ge 2k-1$. 
\end{proof}


\subsection{Grid blocks, Double-forcing blocks, and a proof of Theorem~\ref{klconstruction} (i)}


\subsubsection*{Grid blocks.} 
\begin{figure}
\begin{center}
\epsfig{file={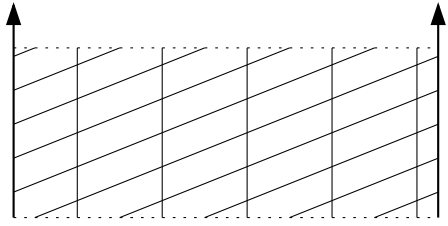}}
\end{center}
\caption{A $(5,1)$-grid block $G(5,1)$.}
\label{figure_block_Cij}
\end{figure}

Let $m, k>0$. 
An {\em $(m,1)$-grid block}, $G(m,1)$ contains two groups, $G'$ and $G''$, 
of cables (or edges). Both $G'$ and $G''$ 
contain $m$ cables.
Refer to
Figure~\ref{figure_block_Cij}. The cables of $G'$ form $m$ parallel
segments in $G(m,1)$. The cables of $G''$ are also parallel in
$G(m,1)$ but make exactly one twist around the cylinder, intersecting every cable of $G'$ 
exactly once. Moreover, the cables from $G'$ and $G''$ intersect 
both the upper and lower boundary alternately.
%
%
An {\em $(m,k)$-grid block} $G(m,k)$ consists of $k$ identical subblocks
$G(m,1)$ stacked on top of each other. 

Observe that $G(2,1)$ is a crossing-forcing configuration and that $G(1, k)$ is a
$k$-spiral. So 
grid blocks are common generalizations 
of the spirals and crossing-forcing blocks.

The following observation generalizes 
Observation~\ref{obs_spiral} and is easily shown by induction. 

\begin{observation}\label{obs_grid}
Every transversal of the grid block $G(m,k)$ has at least $mk-1$ crossings
with the cables in $G(m,k)$. 
\qed
\end{observation}



\subsubsection*{Double-$k$-forcing blocks.} Let $k>0$. A {\em double-$k$-forcing block} $D_k$
contains two groups of cables (edges), say, $D'$ and $D''$. 
\begin{itemize}

\item Each of $D'$ and $D''$ contains $4^k$ cables forming a $k$-forcing block in $B$. 

\item The cables of $D'$ and $D''$ are consecutive along the upper and lower
  boundaries
(but they are ordered differently on the two boundaries). 

\item Any two cables in $D_k$ intersect at most $k$ times.

\end{itemize}

\begin{figure}
\begin{center}
\epsfig{file={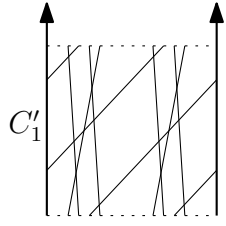}}
\end{center}
\caption{A double-$1$-forcing block.}
\label{figure_hypersubblock}
\end{figure}

We can construct double-$k$-forcing blocks from subblocks in the same way as $k$-forcing blocks.
For $k=1$, the construction is shown on Figure~\ref{figure_hypersubblock}.
Suppose we already have $D_i$. We add a subblock $C'_{i+1}$ to the bottom of $D_i$
as follows.
We divide each cable into four subcables and let these cables form 
a crossing-forcing configuration  
in $C'_{i+1}$, so that any two cables cross at most once in $C'_{i+1}$ and the
subcables of the same cable are consecutive along the upper and lower boundary of $C'_{i+1}$. 

The following is a direct consequence of Lemma~\ref{kforcinglemma}.

\begin{observation}\label{obs_doublekforcing}
Suppose that $D_k$ is a double-$k$-forcing block with 
two groups of cables $D'$ and $D''$. 
Then for any transversal $\kappa$ of $D_k$, there are 
cables $\alpha'\in D'$ and  $\alpha''\in D''$ that both cross $\kappa$ at
least $k$ times. 
\end{observation}

Now, we are ready to prove Theorem~\ref{klconstruction} (i). 
The construction consists of $4k+1$ consecutive blocks, 
$X,\allowbreak D_1,\allowbreak
G_1,\allowbreak D_2,\allowbreak G_2,\allowbreak \ldots,\allowbreak
G_{2k-1},\allowbreak D_{2k},\allowbreak Y$,
in this order (from bottom to top). Let $m=2\cdot 4^k$.  
We define $2km$ 
independent edges
$\alpha_i^j$ and $\beta_i^j$,
$1\le i\le k$, $1\le j\le m$, 
and two isolated vertices $x$ and
$y$ as follows. 
Put $x$ in $X$ and $y$ in $Y$.

\begin{itemize}

\item The edges $\alpha_1^j$ and  $\beta_1^j$, $1\le j\le m$, 
are in $D_1$ and $G_1$. 

\item For every $i$, $1<i<2k$, the edges $\alpha_i^j$ and  
$\beta_i^j$, $1\le j\le m$, are in $G_{i-1}$, $D_i$ and $G_i$.

\item The edges $\alpha_{2k}^j$ and $\beta_{2k}^j$, $1\le j\le m$, are in $G_{2k-1}$ and $D_{2k}$.

\end{itemize}

For every $i$, $1\le i\le 2k$, $D_i$ is a double-$k$-forcing block, and 
the edges $\alpha_i^j$,  $1\le j\le m$, and 
$\beta_i^j$,  $1\le j\le m$, 
form its two groups $D'$ and $D''$.
For every $i$, 
$1\le i\le 2k-1$,
$G_i$ is a $(2,k)$-grid block $G(2,k)$ with groups 
of cables $G'$ and $G''$.
The edges $\alpha_i^j$ form a cable $G'_1$, 
the edges $\beta_i^j$ form a cable  $G'_2$,
the edges $\alpha_{i+1}^j$ form a cable $G''_1$, and  
the edges $\beta_{i+1}^j$ form a cable  $G''_2$.
Cables  $G'_1$ and $G'_2$ form the group $G'$, and cables 
 $G''_1$ and $G''_2$ form the group $G''$. 

Let $\kappa$ be a fixed curve from $x$ to $y$.
We show that $\kappa$ crosses one of the curves 
$\alpha_i^j$ or $\beta_i^j$
at least $2k$ times. 

For every $i$, $1\le i\le 2k$,  by Observation~\ref{obs_doublekforcing},
there is a $j$,
$1\le j\le m$, such that  
$\alpha_i^j$ and $\kappa$ cross at least $k$ times in $D_i$. For simplicity,
denote this 
$\alpha_i^j$ by $\alpha_i$.
Similarly, there is a $j'$,  $1\le j'\le m$, such that  
$\beta_i^{j'}$ and $\kappa$ cross at least $k$ times in $D_i$. Denote
$\beta_i^{j'}$ by $\beta_i$.

For every $Z\in\{D_1,\allowbreak D_2,\allowbreak \ldots ,\allowbreak
D_{2k},\allowbreak G_1,\allowbreak G_2,\allowbreak \ldots ,\allowbreak
G_{2k-1}\}$ and 
$\chi\in\{\alpha_1,\allowbreak \ldots ,\allowbreak \alpha_{2k},\allowbreak
\beta_1,\allowbreak \ldots ,\allowbreak \beta_{2k}\}$,
let $Z(\chi)$ denote the number of intersections of 
$\chi$ with $\kappa$ in $Z$. 
By the choice of $\alpha_i$ and $\beta_i$, for every $i$, $1\le i\le 2k$, we have
$$D_i(\alpha_i), \ D_i(\beta_i)\ge k.$$
By Observation~\ref{obs_grid}, for every $i$, $1\le i\le 2k-1$, we have
$$G_i(\alpha_i)+G_i(\beta_i)+G_i(\alpha_{i+1})+G_i(\beta_{i+1})\ge 2k-1.$$

Summing up, 
\begin{align*}
&\sum_{i=1}^{2k}\left(\cro(\kappa, \alpha_i)+\cro(\kappa, \beta_i)\right)\\
=&\sum_{i=1}^{2k}\left(D_i(\alpha_i)+D_i(\beta_i)\right)+
\sum_{i=1}^{2k-1}\left(G_i(\alpha_i)+G_i(\alpha_{i+1})+G_i(\beta_i)+G_i(\beta_{i+1})\right)\\
\ge &4k^2+(2k-1)^2=4k(2k-1)+1. 
\end{align*}

Therefore, for some $i$, $\cro(\kappa, \alpha_i)\ge 2k$ or $\cro(\kappa,
\beta_i)\ge 2k$. 
\qed


\section{Proof of Theorem~\ref{klconstruction} (ii)}

Let $G$ be a non-complete $k$-simple topological graph,
and let $u$ and $v$ be two non-adjacent vertices of $G$.
We prove that $u$ and $v$ can be connected by a curve that has
at most $2k$ points in common with any edge of $G$.

Place a new vertex at each crossing of $G$ and subdivide the edges accordingly.
Let $G'$ denote the resulting topological (multi)graph. Clearly, there is no loss of generality in assuming that $G'$ is connected. Choose an arbitrary path $\alpha$
in $G'$ connecting $u$ and $v$.
We distinguish two types of vertices on $\alpha$.
A vertex $x$ of $G'$ that lies on $\alpha$ is called a {\em passing vertex\/}
if the two edges of $\alpha$ incident to $x$ belong to the same edge of $G$.
A vertex $x$ of $G'$ that lies on $\alpha$ is a {\em turning vertex\/}
if it is not a passing vertex, that is, if the two edges of $\alpha$ meeting at $x$
belong to distinct edges of $G$.

\begin{figure}[ht]
\begin{center}
\includegraphics{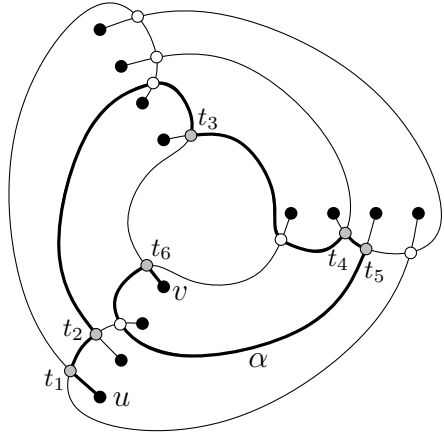}
\caption{A $(u,v)$-path $\alpha$ (in bold) with $c(\alpha)=(6,0,0,1,1,0,1,0)$ and its turning
  vertices $t_i$.\label{codes}}
\end{center}
\end{figure}

Assign to $\alpha$ a unique {\em code}, denoted by $c(\alpha)$, as follows.
Suppose that $\alpha$ contains $r$ turning vertices for some $r\ge 0$.
These vertices divide $\alpha$ into $r+1$ intervals,
$I^{\alpha}_1, I^{\alpha}_2, \ldots ,I^{\alpha}_{r+1}$, ordered from $u$ to $v$.
Set $p^{\alpha}_0=r$ and for any $i$,
$1\le i\le r+1$, let $p^{\alpha}_i$ denote the number of passing vertices on $I^{\alpha}_i$.
Let $c(\alpha)=(p^{\alpha}_0, p^{\alpha}_1, p^{\alpha}_2, \ldots ,
p^{\alpha}_{r+1})$; see Figure~\ref{codes}.

Order the codes of all $(u,v)$-paths lexicographically:
if $\alpha$ and $\beta$ are two $(u,v)$-paths in $G'$,
with codes $c(\alpha)=(p^{\alpha}_0=r, p^{\alpha}_1, p^{\alpha}_2, \ldots , p^{\alpha}_{r+1})$ and
$c(\beta)=(p^{\beta}_0=s, p^{\beta}_1, p^{\beta}_2, \ldots , p^{\beta}_{s+1})$, respectively,
then let $c(\alpha)\prec_{\mbox{\scriptsize{lex}}} c(\beta)$ if and only if
$c(\alpha)\neq c(\beta)$ and for the smallest index $i$ such that $p_i\neq q_i$ we have $p_i<q_i$.

Finally, define a {\em partial ordering} $\prec$ on the set of all the $(u,v)$-paths in $G'$:
for any two $(u,v)$-paths, $\alpha$ and $\beta$, let
$\alpha\prec\beta$ if and only if $c(\alpha)\prec_{\mbox{\scriptsize{lex}}} c(\beta)$.

Let $\gamma$ be a {\em minimal\/} element with respect to $\prec$.
Suppose that $\gamma$ has $r$ turning vertices, $t_1, t_2, \ldots , t_r$,
$r\ge 0$, which divide $\gamma$ into intervals $I^{\gamma}_1, I^{\gamma}_2, \ldots , I^{\gamma}_{r+1}$,
ordered from $u$ to $v$.
Consider the intervals as half-closed, that is, for every $i$, $0\le i \le r$, let $t_i$ belong to $I^{\gamma}_{i+1}$.

Next we establish some simple properties of the intersections of $\gamma$ with the edges of $G$.

\begin{lemma}\label{consecutive}
Let $e$ be an edge of $G$ that has only finitely many points in common with $\gamma$.
Then all of these points belong to two consecutive intervals of $\gamma$.
\end{lemma}

\begin{figure}[h]
\begin{center}
\includegraphics{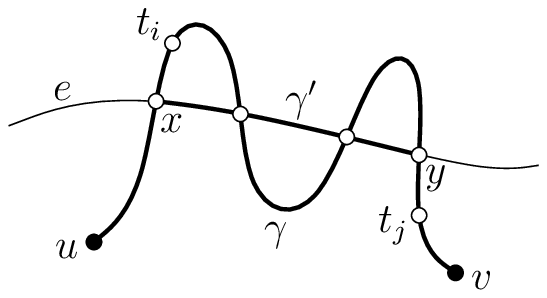}
\caption{Two $(u,v)$-paths $\gamma$ and $\gamma'$ (both in bold) in the proof of Lemma~\ref{consecutive}.\label{modify}}
\end{center}
\end{figure}

\begin{proof}
Suppose for contradiction that $e$ has nonempty intersection with at least
two non-consecutive intervals of $\gamma$.
Let $x$ (and $y$) denote the crossing of $e$ and $\gamma$, closest
to (respectively, farthest from) $u$ along $\gamma$. Let $x$ belong to $I^{\gamma}_i$ and
let $y$ belong to $I^{\gamma}_j$, where $i<j-1$.

Let $\gamma'$ be another $(u,v)$-path, which is identical to $\gamma$
from $u$ to $x$, identical to $e$ from $x$ to $y$, and finally
identical to $\gamma$ from $y$ to $v$; see Figure~\ref{modify}.
If $i<j-2$, then it is evident that $c(\gamma')\prec_{\mbox{\scriptsize{lex}}} c(\gamma)$,
since $\gamma'$ has fewer turning vertices than $\gamma$. If $i=j-2$, then
$\gamma$ and $\gamma'$ have the same number of turning vertices, but
$I^{\gamma'}_i$ contains fewer passing vertices than $I^{\gamma}_i$
(hence $p^{\gamma'}_i < p^{\gamma}_i$), and we have $c(\gamma')\prec_{\mbox{\scriptsize{lex}}} c(\gamma)$.
In both cases we obtain that $\gamma'\prec\gamma$, contradicting the minimality of $\gamma$.
\end{proof}

\begin{lemma}\label{2k-2k-1}
Let $e$ be an edge of $G$ that has only finitely many points in common with $\gamma$.
\begin{itemize}
\item[(i)] If none of the common points is a vertex of $e$, then $e$ crosses $\gamma$ at most $2k$ times.
\item[(ii)] If one of the common points is a vertex of $e$, then $e$ crosses $\gamma$ at most $2k-1$ times.
\end{itemize}
\end{lemma}

\begin{proof}
First, suppose that no vertex of $e$ lies on $\gamma$. By Lemma~\ref{consecutive},
$e$ crosses at most two consecutive intervals of $\gamma$.
Each interval is a part of some edge of $G$ and hence crosses $e$
at most $k$ times. This proves (i).

Suppose next that one of the vertices of $e$ lies on $\gamma$. Observe that such a
vertex must be a turning vertex of $\gamma$, say $t_i$. Again, by Lemma~\ref{consecutive},
$e$ crosses at most two consecutive intervals of $\gamma$. Each interval
is a part of some edge of $G$. Moreover, one of them has a common endpoint
with $e$. Therefore, $e$ crosses one of the intervals at most $k$ times and
the other at most $k-1$ times. This proves (ii).
\end{proof}

Note that no edge $e$ of $G$ that has only finitely many points in common with
$\gamma$ can have both of its endpoints on $\gamma$. Otherwise,
both endpoints must be turning vertices of $\gamma$,
say $t_i$ and $t_j$ for some $i < j$. Since the underlying abstract graph
$G$ is simple (that is, $G$ has no multiple edges), the edge of $G$ that contains
$I^{\gamma}_{i+1}$ must be different from the edge that contains
$I^{\gamma}_j$. Hence, there is at least one turning vertex between
$t_i$ and $t_j$ on $\gamma$. Now consider another $(u,v)$-path $\gamma'$
that is identical to $\gamma$ from $u$ to $t_i$, identical to $e$ from $t_i$ to $t_j$,
and finally identical to $\gamma$ from $t_j$ to $v$. The turning vertices $t_i$ and
$t_{j}$ of $\gamma$ are also turning vertices on $\gamma'$. Since
the turning vertices of $\gamma$ that lie between $t_{i}$ and $t_j$ are
not among the turning vertices of $\gamma'$, $\gamma'$ has fewer turning vertices
than $\gamma$. Therefore, we have $c(\gamma')\prec_{\mbox{\scriptsize{lex}}} c(\gamma)$,
contradicting the minimality of $\gamma$.

\begin{lemma}\label{int}
Let $e$ be an edge of $G$ that contains an interval $I^{\gamma}_i$ of $\gamma$.
Then $e$ and $\gamma$ have at most $k$ points in common outside of $I^{\gamma}_i$.
Furthermore, one of these points is $t_i$, the endpoint of $I^{\gamma}_i$.
\end{lemma}

\begin{proof}
Since $I^{\gamma}_i$ and $I^{\gamma}_{i+1}$ are separated by $t_i$, and
$I^{\gamma}_i$ is contained in $e$, it follows that $e$ cannot contain $I^{\gamma}_{i+1}$.
Similarly, $e$ cannot contain $I^{\gamma}_{i-1}$.

If $e$ has a point $p$ in $I^{\gamma}_j$ with $j<i$, consider another $(u,v)$-path $\gamma'$
that is identical to $\gamma$ from $u$ to $p$, identical to $e$ from $p$ to $t_{i-1}$,
and finally identical to $\gamma$ from $t_{i-1}$ to $v$; see Figure~\ref{int_flips}.
If $j<i-1$, the turning vertices $t_j$ and $t_{i-1}$ of $\gamma$ are not among
the turning vertices of $\gamma'$. Although $p$ was a passing vertex of $\gamma$ and
is now a turning vertex of $\gamma'$, still $\gamma'$ has fewer turning vertices than
$\gamma$. Therefore, $c(\gamma')\prec_{\mbox{\scriptsize{lex}}} c(\gamma)$.
If $j = i-1$, the turning vertex $t_j$ of $\gamma$ is not a turning vertex of
$\gamma'$. Again, $p$ was a passing vertex of $\gamma$ and
is now a turning vertex of $\gamma'$. So, $\gamma$ and $\gamma'$ have the same
number of turning vertices. Since $p$ is not a passing vertex of
$\gamma'$, $I^{\gamma'}_{i-1}$ has fewer passing vertices than $I^{\gamma}_{i-1}$
(hence $p^{\gamma'}_{i-1} < p^{\gamma}_{i-1}$), and we have that $c(\gamma')\prec_{\mbox{\scriptsize{lex}}} c(\gamma)$.
In all of the above cases, we obtain that $\gamma'\prec\gamma$, contradicting the minimality of $\gamma$.

\begin{figure}[h]
\begin{center}
\includegraphics{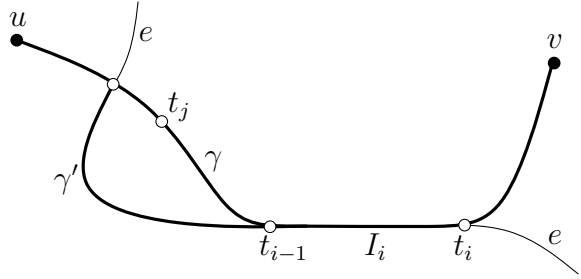}
\caption{Two $u,v$-paths $\gamma$ and $\gamma'$ (both in bold) in the proof of
  Lemma~\ref{int}; $j<i-1$.\label{int_flips}}
\end{center}
\end{figure}

Similarly, if $e$ has a point $p$ in $I^{\gamma}_j$ with $j>i+1$, consider another $(u,v)$-path $\gamma'$
that is identical to $\gamma$ from $u$ to $t_i$, identical to $e$ from $t_i$ to $p$, and finally
identical to $\gamma$ from $p$ to $v$. The turning vertices $t_i$ and $t_{j-1}$ of $\gamma$ are not among
the turning vertices of $\gamma'$. Although $p$ was a passing vertex of $\gamma$ and
is a turning vertex of $\gamma'$, still $\gamma'$ has fewer turning vertices than
$\gamma$. Therefore, $c(\gamma')\prec_{\mbox{\scriptsize{lex}}} c(\gamma)$, contradicting
the minimality of $\gamma$.

Note that the case $j=i+1$ cannot be settled in the same way as the previous cases, since
the number of passing vertices on $e$ between $t_i$ and $p$  may not be smaller than
the number of passing vertices on $\gamma$ between $t_i$ and $p$.
Nevertheless, we can conclude that no interval of $\gamma$ other than $I^{\gamma}_i$ is contained in $e$.
Furthermore, the only interval of $\gamma$ other than $I^{\gamma}_i$ that can share some points
with $e$ is $I^{\gamma}_{i+1}$. Let $f$ be the edge of $G$ that contains $I^{\gamma}_{i+1}$.
Since $e$ and $f$ have at most $k$ points in common, $e$ and $I^{\gamma}_{i+1}$ can have
at most $k$ points in common, too. The point $t_{i}$, the common endpoint of $I^{\gamma}_i$
and $I^{\gamma}_{i+1}$, is one of these points.
\end{proof}

Now we are in a position to complete the proof of 
Theorem~\ref{klconstruction} (ii).
Join $u$ and $v$ by a curve $\beta$ that runs very close to $\gamma$.

We claim that any edge $e$ of $G$ has at most $2k$ points in common
with $\beta$. If $e$ has only finitely many points in common with $\gamma$
and none of them is a vertex of $e$, then every crossing between $e$ and $\beta$
corresponds to a crossing between $e$ and $\gamma$. Therefore, by Lemma~\ref{2k-2k-1}(i),
$e$ and $\beta$ cross each other at most $2k$ times.
If $e$ has only finitely many points in common with $\gamma$, but
one of them is a vertex of $e$, then each crossing between
$e$ and $\beta$ corresponds to a crossing between $e$ and $\gamma$, and there
may be an additional crossing near the vertex of $e$ on $\gamma$.
Again, by Lemma~\ref{2k-2k-1}(ii), there are at most $2k$ crossings between $e$
and $\beta$. Finally, if $e$ contains a whole interval $I^{\gamma}_i$ of
$\gamma$, then each crossing between $e$ and $\beta$
corresponds to a crossing between $e$ and $\gamma$, or to a vertex of
$e$ on $\gamma$. There may be an additional crossing near the endpoint
$t_{i}$ of $I^{\gamma}_i$. Thus, there are at most $k+1$ crossings.
\qed


\section{Proof of Theorem~\ref{linearbound}: Upper Bounds} 

The construction for $k=1$ is essentially 
different from the constructions for
$k>1$. 
For $k>1$, all constructions are variations of the constructions used in the
proofs of Theorem~\ref{klconstruction}, Lemma~\ref{lemma_weakklconstruction} and
Lemma \ref{2k-1construction}, but they give different bounds for different values of
$k$.
Table~\ref{table1} shows our best upper bounds for different values of $k$.

\begin{figure}[h]
\begin{center}
\includegraphics{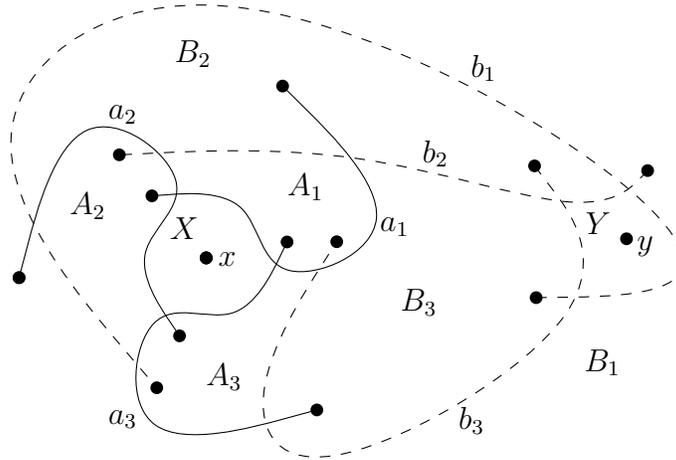}
\caption{A topological graph $G$: the edge $\{x, y\}$ cannot be added.\label{nembovitheto1}}
\end{center}
\end{figure}

\paragraph{First construction.}
This construction is for $k=1$.
First we need to modify the graph $G_1$ on Figure~\ref{kampok}.
Consider the edges of $G_1$ incident to $x$, and modify them
in a small neighborhood of $x$ so that the resulting edges have distinct
endpoints, they pairwise cross each other, and their union encloses a {\em region} $X$
(i.e., a connected component $X$ of the complement of the union of the edges) which
contains $x$. Analogously, modify the other three
edges of $G_1$ in a small neighborhood of $y$. Let $Y$ be the region that
contains $y$ and is enclosed by the modified edges. The resulting
simple topological graph $G$ has $12$ vertices and $6$ edges; see
Figure~\ref{nembovitheto1}. The points $x,y\in V(G_1)$ do not belong to $V(G)$.

\begin{lemma}\label{2cr}
Let $x$ and $y$ be any pair of points belonging to the regions $X$ and $Y$ in $G$,
respectively. Then any curve joining $x$ and $y$ will meet at least one of the edges of $G$ at
least twice.
\end{lemma}

\begin{proof}
We prove the claim by contradiction.
Let $a_1, a_2, a_3$, $b_1, b_2, b_3$ denote the
edges of $G$. They divide the plane into eight regions, $X$, $Y$, $A_1, A_2,
A_3$, $B_1, B_2, B_3$; see Figure~\ref{nembovitheto1}.
Suppose there exists an oriented curve from $x$ to $y$
that crosses every edge of $G$ at most once. Let $\gamma$ be such a curve with
the smallest number of crossings with the edges of $G$.
Let $c_1, c_2, \ldots , c_{m-1}$ be the crossings between $\gamma$ and the edges of
$G$, ordered according to the orientation of $\gamma$. They divide
$\gamma$ into intervals $I_1, I_2, \ldots , I_m$, ordered again according to the orientation of
$\gamma$. The first interval $I_1$ lies in $X$, and the last one, $I_m$, lies in $Y$.
Observe that no other interval can belong to $X$ or to $Y$, because in this case we
could simplify $\gamma$ and obtain a curve with a smaller number of crossings.
By symmetry, we can assume that the first crossing, $c_1$, is a crossing
between $\gamma$ and $a_1$. Then $I_2$ belongs to $A_1$. The following property holds.

\medskip

{\sl Property $\P$: If for some $j\ge 2$, the interval $I_j$ belongs to $A_i$
  (or $B_i$), then one of the points $c_1, c_2, \ldots , c_{j-1}$ is a crossing
  between $\gamma$ and the edge $a_i$ (or $b_i$, respectively).}

\medskip

We prove Property $\P$ by induction on $j$. Clearly, the property holds for
$j = 2$. Assume that $I_{j-1}$ is in $A_i$ (or $B_i$)
and one of $c_1, c_2, \ldots , c_{j-2}$ is a crossing between $\gamma$ and $a_i$
(or $b_i$). For simplicity, assume that $I_{j-1}$ belongs to the region $A_1$
and that one of the points $c_1, c_2, \ldots , c_{j-2}$ is a crossing between $\gamma$ and
$a_1$; the other cases are analogous. Since $c_{j-1}$ cannot belong to
$a_1$, it must be a crossing between $\gamma$ and either $a_2$ or $b_2$. In the
first case, $I_j$ belongs to $A_2$, in the second to $B_2$. In either case, Property $\P$
is preserved.

Now, we can complete the proof of Lemma~\ref{2cr}.
Consider the interval $I_{m-1}$. Since $I_m$ lies in $Y$, for some $i$, the interval
$I_{m-1}$ must lie in $B_i$. Suppose for simplicity that
$I_{m-1}$ lies in $B_1$. By Property $\P$ (with $j = m-1$, $m\ge 3$),
one of the points $c_1, c_2, \ldots , c_{m-2}$ must be a crossing between $\gamma$ and $b_1$.
However, using that $I_m$ is in $Y$, $c_{m-1}$ must be another crossing between $\gamma$ and $b_1$.
Thus, $\gamma$ crosses $b_1$ twice, which is a contradiction.
\end{proof}

\begin{figure}[h]
\begin{center}
\includegraphics{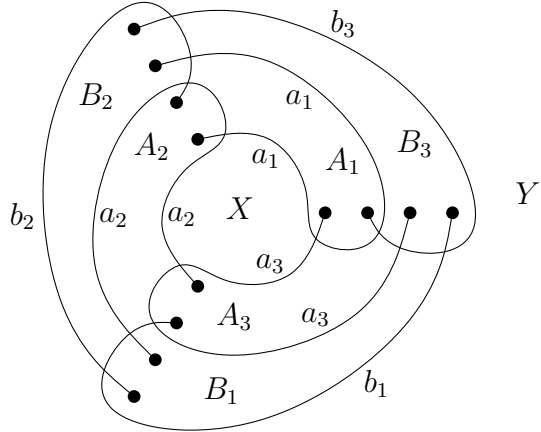}
\caption{A topological graph $H$, a modification of $G$.\label{nembovitheto-gyurus}}
\end{center}
\end{figure}

Now, we return to the proof of the upper bound in Theorem~\ref{linearbound}.
Modify the drawing of $G$ in Figure~\ref{nembovitheto1} so that the region $Y$
becomes unbounded, and let $H$ be the resulting topological graph.
Denote by $Y$ the {\em outer region\/} of $H$ and by $X$ the {\em inner region\/} of $H$; see Figure~\ref{nembovitheto-gyurus}.

\begin{figure}[h]
\begin{center}
\scalebox{0.4}{\includegraphics{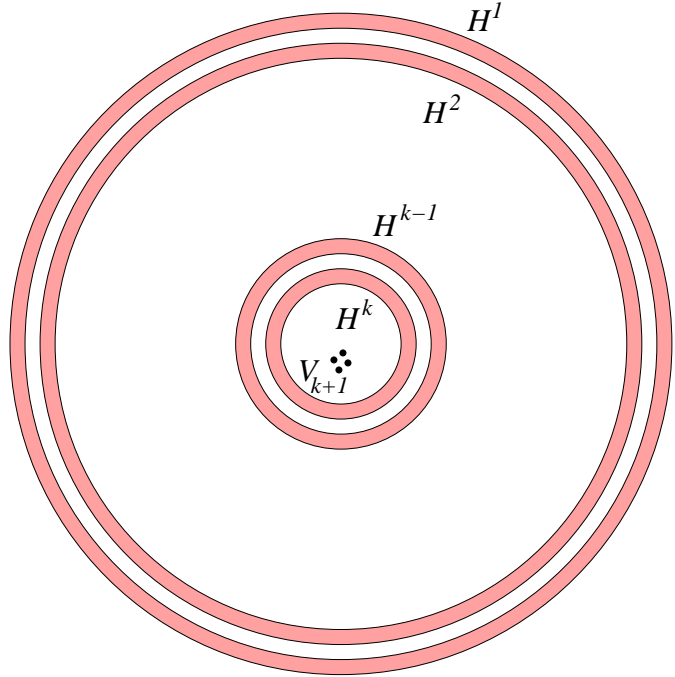}}
\caption{A saturated simple topological graph $F_n$. \label{soknembovitheto}}
\end{center}
\end{figure}

For every $n\ge 1$, construct a saturated simple topological graph
$F_n$, as follows. Let $k=\lfloor{n/12}\rfloor$.
Take a disjoint union of $k$ scaled and translated copies
of $H$, denoted by $H^1,\allowbreak H^2,\allowbreak \ldots,\allowbreak H^k$, such that for any $i$, $1<i\le k$,
the copy $H^i$ lies entirely in the inner region of $H^{i-1}$; see Figure~\ref{soknembovitheto}.
For $1\le i\le k$, let $V_i$ be the vertex set of $H^i$. Finally, place
$n-12k$ additional vertices in the inner region of $H^k$, and let $V_{k+1}$
denote the set of these vertices.
Obviously, we have $|V_{k+1}| < 12$.

Add to this topological graph all possible missing edges one by one, in an
arbitrary order, as long as
it remains simple. 
We end up with a saturated simple topological graph $F_n$
with $n$ vertices.
Observe that for every $i$ and $j$ with $1\le i<j-1<k$, $V_i$ lies in the outer region
of $H^{i+1}$, while $V_j$ is in the inner region of $H^{i+1}$.
By Lemma~\ref{2cr} (applied with $G = H^{i+1}$, $x\in V_j$, $y\in V_i$),
no edge of $F_n$ runs between $V_i$ and $V_j$.
Hence, every vertex in $V_i$ can be adjacent to at most $35$ other vertices; namely,
to the elements of $V_{i-1}\cup V_i\cup V_{i+1}$. Therefore,
$F_n$ is a saturated simple topological graph with $n$ vertices and at most
$17.5n$ edges.

\begin{figure}
\begin{center}
\epsfig{file={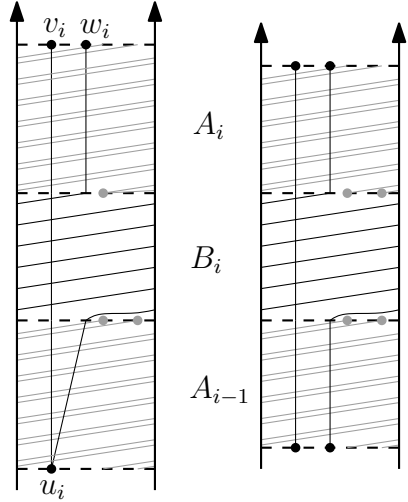}}
\end{center}
\caption{Three consecutive blocks in the constructions of saturated $k$-simple
topological graphs. Left: the second construction for $k=7$. Right: the third 
construction for $k=6$.}
\label{figure_saturated_first_second}
\end{figure}

\paragraph{Second construction.}
This construction is used for all odd $k\ge 5$ and all even $k\ge 12$.
Suppose for simplicity that $n$ is divisible by $3$ and let
$m=n/3$.  The construction consists of $2m+3$ consecutive blocks, 
$B_0,\allowbreak A_0,\allowbreak B_1,\allowbreak
A_1,\allowbreak \ldots,\allowbreak B_{m},\allowbreak A_{m},\allowbreak B_{m+1}$, 
in this order, from bottom to top. See
Figure~\ref{figure_saturated_first_second}, left.

For every $i$, $1\le i\le m$, let $u_i$ be a vertex on the common boundary of
$B_{i-1}$ and 
$A_{i-1}$ and let $v_i$ and $w_i$ be vertices on the common boundary of
$A_{i}$ and $B_{i+1}$.
Let $\alpha_i$ be an edge connecting $u_i$ and $v_i$ and let $\beta_i$ be an
edge connecting $u_i$ and
$w_i$. The pair $(\alpha_i, \beta_i)$ is called the {\em $i$th bundle}. 
The edges $\alpha_i$ and $\beta_i$ form a $(k-1)$-spiral in $B_i$.
For $1< i\le m$, the edges $\alpha_i$ and $\beta_i$ form a cable in $A_{i-1}$, the
edges $\alpha_{i-1}$ and $\beta_{i-1}$ also form a cable in $A_{i-1}$, and
these two cables form a $k$-spiral in $A_{i-1}$. 
The resulting $k$-simple topological graph $G$ has $n$ vertices and
$2n/3$ edges.
Add to $G$ all possible missing edges one by one, as long
as the drawing remains $k$-simple.
We obtain a saturated $k$-simple topological graph $H$. Note that $H$ is not
uniquely determined by $G$, not even as an abstract graph.

Suppose that $k\ge 11$. Just like in the proof of 
Lemma~\ref{lemma_weakklconstruction}, 
we can prove that any curve from $A_i$ to $A_{i+3}$ has to cross one
of the curves
$\alpha_{i+1}$, $\beta_{i+1}$, $\alpha_{i+2}$, $\beta_{i+2}$, $\alpha_{i+3}$,
$\beta_{i+3}$
at least $(k-2)/2+2(k-1)/3>k$ times. 
Therefore, in $H$, a vertex from the block $A_i$ can possibly be connected
only to other vertices from the five blocks $A_{i-2},\allowbreak
A_{i-1},\allowbreak A_i,\allowbreak A_{i+1},\allowbreak A_{i+2}$. Since every
block $A_j$ has at most three vertices, the maximum degree in $H$ is at most
$5\cdot3-1=14$ and thus $H$ has at most $7n$ edges. 

For $k \ge 9$ odd and for every $j$, $1\le j\le 3$, 
every curve $\kappa$ from $A_i$ to $A_{i+3}$ has to cross one
of the two curves $\alpha_{i+j},\allowbreak \beta_{i+j}$ at least $(k-1)/2$
times in $B_{i+j}$. Let $\gamma_{i+j}$ be this curve. Now, for every $j$, $1\le
j\le 2$, the curve $\kappa$ crosses $\gamma_{i+j}$ and $\gamma_{i+j+1}$
together at least $k-1$ times in $A_{i+j}$. It follows that $\kappa$ crosses
one of the curves $\gamma_{i+1},\allowbreak \gamma_{i+2},\allowbreak
\gamma_{i+3}$ at least $(1/2 + 2/3)\cdot (k-1)>k$ times. 

Therefore, in $H$, a vertex from the block $A_i$ can be connected only to
other vertices from the five blocks $A_{i-2},\allowbreak A_{i-1},\allowbreak
\dots,\allowbreak A_{i+2}$. The maximum degree in $H$ is thus at most $5\cdot
3 - 1 =14$ and $H$ has at most $7n$ edges.

Similarly, for $k \ge 7$ odd, for every curve $\kappa$ from $A_i$ to $A_{i+4}$,
there are four curves $\gamma_{i+1},\allowbreak \gamma_{i+2},\allowbreak
\gamma_{i+3},\allowbreak \gamma_{i+4}$ such that $\kappa$ crosses one of them
at least $(1/2 + 3/4)\cdot (k-1)>k$ times. 

Therefore, in $H$, a vertex from the block $A_i$ can be connected only to
vertices from the seven blocks $A_{i-3},\allowbreak A_{i-2},\allowbreak
\dots,\allowbreak A_{i+3}$. The maximum degree in $H$ is thus at most $7\cdot
3 - 1 =20$ and $H$ has at most $10n$ edges.

For $k \ge 5$ odd, for every curve $\kappa$ from $A_i$ to $A_{i+5}$, there are
five curves, $\gamma_{i+1},\allowbreak \gamma_{i+2},\allowbreak
\dots,\allowbreak \gamma_{i+5}$ such that $\kappa$ crosses one of them at
least $(1/2 + 4/5)\cdot (k-1)>k$ times. 

Therefore, in $H$, a vertex from the block $A_i$ can be connected only to
vertices from the nine blocks $A_{i-4},\allowbreak A_{i-3},\allowbreak
\dots,\allowbreak A_{i+4}$. The maximum degree in $H$ is thus at most $9\cdot
3 - 1 =26$ and $H$ has at most $13n$ edges.

\paragraph{Third construction.}
This construction is used for $k\in\{4,\allowbreak 6,\allowbreak 8,\allowbreak
10\}$. It is a modification of the second construction, where the edges of the
$i$th bundle, 
$\alpha_i$ and $\beta_i$, do not have common endpoints, so they form a matching
rather then a path, and they form a $k$-spiral in $B_i$. See
Figure~\ref{figure_saturated_first_second}, right.

For $k \ge 6$ even and for every $j$, $1\le j\le 3$, 
every curve $\kappa$ from $A_i$ to $A_{i+3}$ has to cross one
of the two curves $\alpha_{i+j},\allowbreak \beta_{i+j}$ at least $k/2$ times
in $B_{i+j}$. Let $\gamma_{i+j}$ be this curve. Now, for every $j$, $1\le j\le
2$, the curve $\kappa$ crosses $\gamma_{i+j}$ and $\gamma_{i+j+1}$ together at
least $k-1$ times in $A_{i+j}$. It follows that $\kappa$ crosses one of the
curves $\gamma_{i+1},\allowbreak \gamma_{i+2},\allowbreak \gamma_{i+3}$ at
least $k/2 + 2(k-1)/3>k$ times. 

Therefore, in $H$, a vertex from the block $A_i$ can be connected only to
other vertices from the five blocks $A_{i-2},\allowbreak A_{i-1},\allowbreak
\dots,\allowbreak A_{i+2}$. Since every block $A_j$ now has at most four
vertices, the maximum degree in $H$ is at most $5\cdot 4 - 1 =19$ and $H$ has
at most $9.5n$ edges.

Similarly for $k \ge 4$ even, for every curve $\kappa$ from $A_i$ to
$A_{i+4}$, there are four curves $\gamma_{i+1},\allowbreak
\gamma_{i+2},\allowbreak \gamma_{i+3},\allowbreak \gamma_{i+4}$ such that
$\kappa$ crosses one of them at least $k/2 + 3(k-1)/4>k$ times. 

Therefore, in $H$, a vertex from the block $A_i$ can be connected only to
other vertices from the seven blocks $A_{i-3},\allowbreak A_{i-2},\allowbreak
\dots,\allowbreak A_{i+3}$. The maximum degree in $H$ is thus at most $7\cdot
4 - 1 =27$ and $H$ has at most $13.5n$ edges.

\paragraph{Fourth construction.}
This construction is for $k=2$. 
First we present a weaker but simpler version.
It is a modification of the previous constructions.
Here each bundle contains $16$ independent edges.
The edges of the $i$th bundle form a $2$-forcing block
in  $B_i$.
In $A_i$, the edges of the $i$th bundle form a cable, the edges 
of the $(i+1)$st bundle form another cable, and these two cables form 
a $2$-spiral. 
Let $\kappa$ be a curve from $A_i$ to $A_{i+2}$.
Just like in the previous arguments, using Observation~\ref{obs_spiral} and 
Lemma~\ref{kforcinglemma}, 
it is not hard to see that 
$\kappa$ has to cross an edge more than twice.
Therefore, in $H$, a vertex from $A_i$ can be connected only to
vertices from $A_{i-1},\allowbreak A_{i},\allowbreak A_{i+1}$. 
Every block $A_j$ has at most $32$
vertices, so the maximum degree in $H$ is at most $95$, therefore, $H$ has
at most $47.5n$ edges.  

\begin{figure}
\begin{center}
\epsfig{file={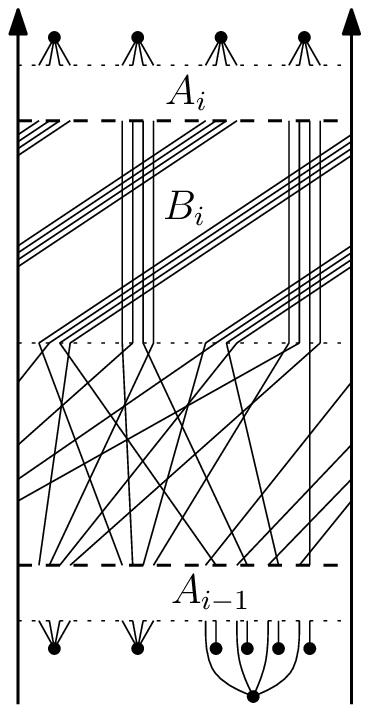}} \hskip 1cm
\epsfig{file={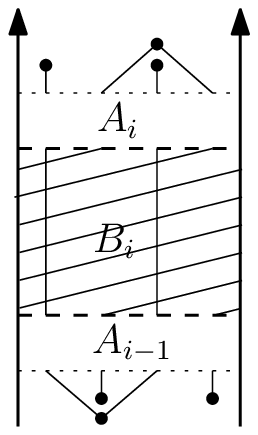}}
\end{center}
\caption{A bundle in the construction of a saturated $2$-simple (left) and a
saturated $3$-simple (right) topological graph.}
\label{figure_saturated_for_k_4_2}
\end{figure}

The best construction we have is very similar. To obtain it, in each bundle we identify
some of the endpoints of the edges, and we also modify the order of the edges along the bottom boundary of $B_i$; see Figure~\ref{figure_saturated_for_k_4_2}, left.
Then every block $A_i$ has at most $11$
vertices, so the maximum degree in $H$ is at most $3\cdot 11-1=32$ and $H$ has at most $16n$ edges.

\paragraph{Fifth construction.}
This construction is for $k=3$.  
First we present a weaker but simpler version.
It is again a modification of the previous constructions.
Here each bundle contains four independent edges.
The edges of the $i$th bundle form a grid block $G(2,3)$ in  $B_i$.
In $A_i$, the edges of the $i$th bundle form a cable, the edges 
of the $(i+1)$st bundle form another cable, and these two cables form 
a $3$-spiral. 
Let $\kappa$ be a curve from $A_i$ to $A_{i+3}$.
Just like in the previous arguments, using Observations~\ref{obs_grid} and~\ref{obs_spiral}, 
it is not hard to see that 
$\kappa$ has to cross an edge more than three times.
Therefore, in $H$, a vertex from $A_i$ can be connected only to
vertices from $A_{i-2},\allowbreak A_{i-1},\allowbreak
\dots,\allowbreak A_{i+2}$. Every block $A_j$ has at most $8$
vertices, so $H$ has
at most $19.5n$ edges.  

To obtain our best construction, in each bundle we identify
some endpoints
of the edges; see  Figure~\ref{figure_saturated_for_k_4_2}, right.
Then every block $A_i$ has at most $6$
vertices, so the maximum degree in $H$ is at most $5\cdot 6-1=29$ and  $H$
has at most
$14.5n$ edges.
A modification of this construction works for $k=1$, and it gives 
the same upper bound, $17.5n$, as the first construction. 

\medskip

This concludes the proof of the upper bounds.
\qed

\section{Proof of Theorem~\ref{linearbound}: Lower Bounds}

A vertex of a (topological) graph is
{\em isolated\/} if its degree is zero. A triangle in a (topological) graph is
called {\em isolated\/} if its vertices are incident to no edges other
than the edges of the triangle.

\begin{lemma}\label{triangle}
A saturated simple topological graph on at least four vertices contains no isolated triangle.
\end{lemma}

\begin{figure}[h]
\begin{center}
\includegraphics{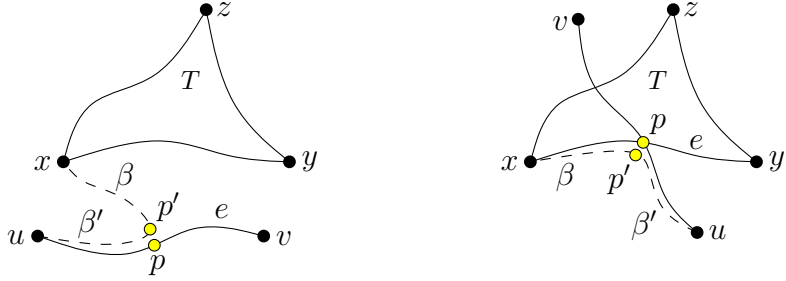}
\caption{Case 1 and Case 2 of Lemma~\ref{triangle}.\label{trianglecases}}
\end{center}
\end{figure}

\begin{proof}
Let $G$ be a saturated simple topological graph with at least four vertices, and suppose
for contradiction that $G$ has an
isolated triangle $T$ with vertices $x$, $y$ and $z$. By definition, the edges of $T$
do not cross one another. 

If all vertices other than $x,y,z$ are isolated, it is trivial to add a new edge without crossings. Hence we may assume that $G$ has an edge not contained in $T$.
We distinguish two cases.

\smallskip

{\bf Case 1.} The edges of $T$ cross no other edges.

The edges of $G$ divide the plane into regions. Let $R$ denote a region
bounded by the edges of $T$ and at least one other nontrivial curve $\omega$.
Let $e = \{u,v\}$ be an edge that contributes to $\omega$, and let $p$ be a point on $e$
that belongs to the boundary of $R$; see Figure~\ref{trianglecases}, left.
Choose a point $p'$ inside of $R$, very close to $p$.
Let $\beta$ be a curve running inside $R$ that connects a vertex of $T$, say $x$, to $p'$.
Let $\beta'$ be a curve joining $p'$ and $u$, and running very close to the edge $e$.
Adjoining $\beta$ and $\beta'$ at $p'$, we obtain a curve $\gamma$ connecting $x$ and $u$, two previously non-adjacent vertices of $G$.
The curve $\gamma$ crosses neither an edge of $T$ or an edge of $G$ incident to $u$.
Since $\beta$ is crossing-free, all crossings between $\gamma$ and the edges of $G$
must lie on $\beta'$ and, hence, must correspond to crossings along the
edge $e$. Therefore, every edge of $G$ can cross $\gamma$ at most once.
Consequently, $\gamma$ can be added to $G$ as an extra edge so that the topological graph remains simple.
This contradicts the assumption that $G$ was saturated.

\smallskip

{\bf Case 2.} At least one edge of $T$ participates in a crossing.

Assume without loss of generality that $e=\{x,y\}$ is crossed by another
edge of $G$. Let $p$ denote the crossing on $e$ {\em closest\/} to $x$, and
suppose that $p$ is a crossing between $e$ and another edge $f=\{u,v\}$; see Figure~\ref{trianglecases}, right.
The point $p$ divides $f$ into two parts. At least
one of them, say, $up$, does not cross the edge $\{x,z\}$ of $T$.
The edges $e$ and $f$ divide a small neighborhood of $p$ into four parts.
Choose a point $p'$ in the part bounded by $up$ and $xp$.
Let $\beta$ be a curve connecting $x$ and $p'$, running very close to $e$.
Let $\beta'$ be a curve between $p'$ and $u$, running very close to $f$.
Adjoining $\beta$ and $\beta'$ at $p'$ we obtain a curve $\gamma$
connecting $x$ and $u$, two vertices that were not adjacent in $G$.
Just like in the previous case, add $\gamma$ to $G$ as an extra edge.
The curve $\gamma$ crosses no edge incident to $x$ or $u$.
Since the portion $xp$ of $e$ is crossing-free, $\beta$ must be crossing-free, too.
Therefore, all possible crossings between $\gamma$ and the edges of $G$
must lie on $\beta'$ and, hence, correspond to crossings along $f$.
Thus, every edge of $G$ crosses $\gamma$ at most once, contradicting our
assumption that $G$ was saturated.
\end{proof}

\begin{lemma}\label{isol+deg1}
For any $k>0$, a saturated $k$-simple topological graph on at least three vertices contains
\item[(i)] no isolated vertex,
\item[(ii)] no vertex of degree one.
\end{lemma}


The proof of Lemma~\ref{isol+deg1} is very similar to the proof of Lemma~\ref{triangle}, but much easier. We omit the details.

\smallskip

The lower bound in Theorem~\ref{linearbound} (ii) now follows directly. In a
saturated 
$k$-simple topological graph on $n$ vertices, every vertex has degree at least two, therefore,
it has at least $n$ edges. 
We are left with the proof of the lower bound of part (i). 
It follows immediately 
from the statement below.

\begin{lemma}\label{deg3}
In every saturated simple topological graph with at least four vertices, every vertex has degree at least 3.
\end{lemma}

\begin{proof}
We prove the claim by contradiction. Let $G$ be a saturated simple topological graph,
and let $x$ be a vertex of degree two in $G$. (By Lemma~\ref{isol+deg1}, the degree of
$x$ cannot be $0$ or $1$.)
Let $y$ and $z$ denote the neighbors of $x$. By definition, the edges $\{x,y\}$ and $\{x, z\}$ do not cross.
We distinguish two cases.

\smallskip

{\bf Case 1.} The edges $\{x,y\}$ and $\{x, z\}$ cross no other edges.

By Lemma~\ref{triangle} and~\ref{isol+deg1}, $y$ and $z$ both have degree at least two,
and $x$, $y$ and $z$  do not span an isolated triangle. Hence,
at least one of the vertices $y$ and $z$, say, $y$, has a neighbor $w$
different from $x$ and $z$.
Let $\gamma$ be a curve connecting $x$ to $w$ that runs very close to the
edge $\{x,y\}$ from $x$ to a point in a small neighborhood of $y$, and
from that point all the way to $w$ very close to the edge $\{y,w\}$.
We can assume that $\gamma$ does not cross $\{x,y\}$ and $\{y,w\}$.
Add $\gamma$ to $G$ as an extra edge. Clearly, $\gamma$ crosses no edge
incident to $x$ or $w$, and crosses no edge of $G$ twice.
This contradicts the assumption that $G$ was saturated.

\smallskip

{\bf Case 2.} At least one of the edges $\{x,y\}$ and $\{x, z\}$ participates in a crossing.

Assume without loss of generality that $e = \{x,y\}$ is crossed by another
edge of $G$. Let $p$ be the crossing on $e$ {\em closest\/} to $x$, and
suppose that the other edge passing through $p$ is $f=\{u,v\}$.
The point $p$ divides $f$ into two pieces, at least
one of which, say, $up$, has no point in common with the edge $\{x,z\}$.
Let $\gamma$ be a curve connecting $x$ and $u$, following $e$ very
closely from $x$ to a point in a small neighborhood of $p$, and
from that point following $f$ all the way to $u$.
We can assume that $\gamma$ does not cross $e$ and $f$.
Add $\gamma$ to $G$ as an extra edge. It is again easy to see that this new edge meets no original edge of $G$ more than once, and again, this contradicts the assumption that $G$ was
saturated.
\end{proof}


\section{Proof of Theorem~\ref{3k/2saturated}}

We start with a piece of the construction we used in the 
proof of Theorem~\ref{klconstruction} (i). Then we add some edges
so that it remains a $k$-simple topological graph, and we show that 
it is $(k, \lceil 3k/2\rceil)$-saturated.

Let $D_1$, $G$, $D_2$ be three consecutive blocks, say, from bottom to top,
and let $m=4^k$.
We define $4m$ independent edges 
$\alpha_i^j$, $\beta_i^j$,
$1\le i\le 2$, $1\le j\le m$.
The edges $\alpha_1^j$ and  $\beta_1^j$, $1\le j\le m$, 
are in $D_1$ and $G$, with endpoints on the lower boundary of 
$D_1$ and the upper boundary of $G$.
Denote the sets of these vertices by $V_0$ and $V_2$, respectively.
The edges $\alpha_2^j$ and  $\beta_2^j$, $1\le j\le m$,
are in $G$ and $D_2$, with endpoints on the lower boundary of 
$G$ and the upper boundary of $D_2$.
Denote the sets of these vertices by $V_1$ and $V_3$, respectively.

For $i=1, 2$, 
the block $D_i$ is a double-$k$-forcing block, 
the edges $\alpha_i^j$,  $1\le j\le m$, and 
$\beta_i^j$,  $1\le j\le m$, 
form its two groups $D'$ and $D''$.
The block $G$ is a $(2,k)$-grid block $G(2,k)$ with groups 
of cables $G'$ and $G''$.
The edges $\alpha_1^j$ form a cable $G'_1$, 
the edges $\beta_1^j$ form a cable  $G'_2$,
the edges  $\alpha_{2}^j$ form a cable $G''_1$, and  
the edges  $\beta_{2}^j$ form a cable  $G''_2$.
The cables  $G'_1$ and $G'_2$ form the group $G'$, and the cables 
$G''_1$ and $G''_2$ form the group $G''$ in $G$.
Let $T$ denote the resulting topological graph. 

Let $v_0\in V_0$ and $v_3\in V_3$ be arbitrary vertices, and let 
$\kappa$ be a curve connecting $v_0$ and $v_3$. 
By Observation~\ref{obs_doublekforcing} 
there are edges $\alpha_1=\alpha_1^j$ and $\beta_1=\beta_1^{j'}$
that both cross $\kappa$ at least $k$ times in $D_1$.
Similarly, there are edges $\alpha_2=\alpha_2^l$ and $\beta_2=\beta_2^{l'}$
that both cross $\kappa$ at least $k$ times in $D_2$.
Since $\alpha_1, \allowbreak \alpha_2, \allowbreak
\beta_1, \allowbreak \beta_2$ form a $(2,k)$-grid block in $G$,
by Observation~\ref{obs_grid}, $\kappa$ crosses them in $G$ together 
at least $2k-1$ times. 
Therefore, 
$\kappa$ crosses one of the curves at least 
$\lceil(6k-1)/4\rceil =\lceil 3k/2\rceil$ times.

\begin{figure}
\begin{center}
\epsfig{file={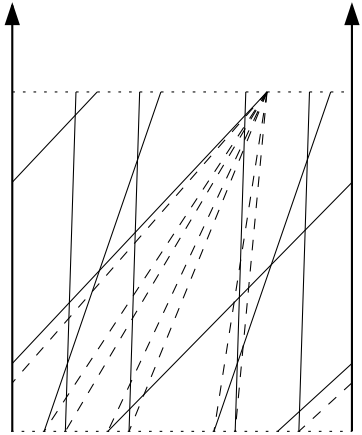}}
\end{center}
\caption{Adding edges to the subblock $C_k$.}
\label{figure_saturated_jump}
\end{figure}

Now we show that any two vertices $v_i\in V_i$ and $v_i\in V_j$ with $|i-j|\le 2$, 
can be connected so that we still have a $k$-simple topological graph.
We only sketch the argument.
By definition, block $D_1$ is divided into $k$ subblocks, $C'_1, C'_2,
\ldots, C'_k$, from top to bottom. Let $v_0\in V_0$, $v_2\in V_2$, 
let $\alpha$ be the edge of $T$ incident with $v_0$ and let 
$\beta$ be the edge of $T$ incident with $v_2$. We can assume that $\alpha\neq\beta$, otherwise we
are done.
Draw a curve $\kappa$ from $v_2$ very close to $\beta$ all the way in $G$, 
and then in the subblocks $C'_1, C'_2,\ldots, C'_{k-1}$. 
In the last subblock, $C'_k$, connect $\kappa$ to $v_0$ so that it crosses 
all edges at most once in $C'_k$. A straightforward but slightly 
technical 
argument shows that 
it is possible. 
For example, we can draw the cables in $C'_k$ as in
Figure~\ref{figure_saturated_jump}, and
then draw $\kappa$ as the shortest line with positive slope.
Repeat this procedure for all pairs $v_0\in V_0$, $v_2\in V_2$. The resulting
topological graph is still $k$-simple. 
We can add similarly all edges between vertices $v_1\in V_1$ and $v_3\in V_3$.
We can connect all the remaining pairs, $v_i\in V_i$ and $v_j\in V_j$, $|i-j|\le 1$,
in a similar, but simpler way.  
We obtain a $(k, \lceil 3k/2\rceil)$-saturated topological graph. \qed

\section{Proof of Theorem~\ref{pseudoline}}

Let $\alpha$ and $\beta$ be two simple closed curves in the plane.
Suppose that $\alpha$ contains $\beta$ in its interior. 
The region between $\alpha$ and $\beta$ is called an 
{\em annulus}. It is homeomorphic to the cylindrical surface, 
so we can transform  {\em blocks\/} onto the annulus.
The region outside $\alpha$ is called the {\em outer
exterior\/} of the annulus. Similarly, the region inside $\beta$ is called the {\em inner exterior\/} of the annulus.

It is enough to prove the following statement; Theorem~\ref{pseudoline}
easily follows.

\begin{theorem}\label{19}
\begin{enumerate}
\item[(i)] For $m\ge 2$ and $k=4m$, there is a $(3,5k/4-5)$-forcing
  $k$-pseudoline arrangement. 
\item[(ii)] For $m\ge 3$ and $k=2^m$, there is a 
$(k,(k/2-2)\cdot (\log_2k +1))$-forcing $k$-pseudoline arrangement. 
\end{enumerate}
\end{theorem}

\begin{figure}
\begin{center}
\epsfig{file={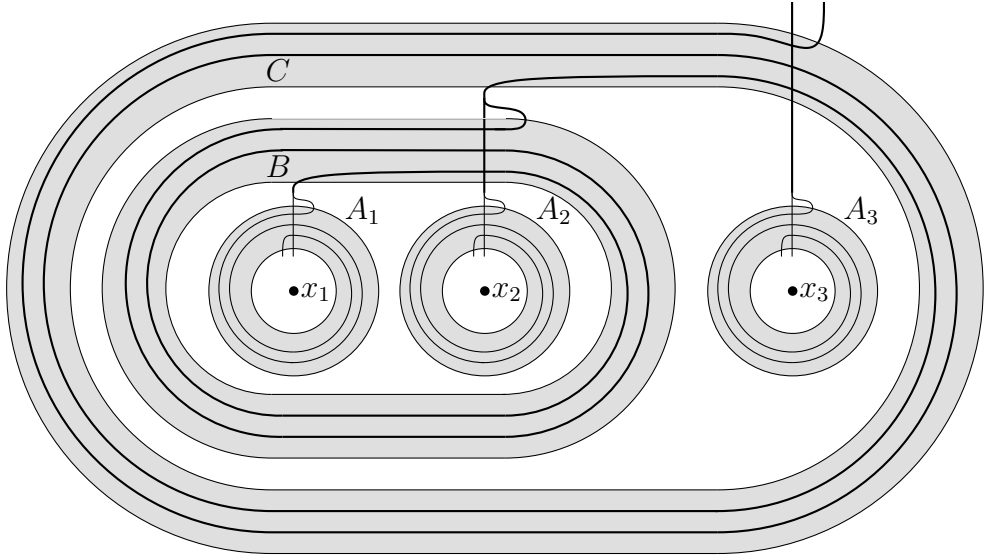}}
\end{center}
\caption{The arrangement from the proof of Theorem~\ref{19} (i).}
\label{figure_k_pseudorays}
\end{figure}

\begin{proof}
(i) Refer to Figure~\ref{figure_k_pseudorays}. 
First 
we construct an arrangement of one-way infinite curves. 
Let $x_1$, $x_2$ and $x_3$ be three distinct points in the plane. 
Let $A_1$, $A_2$ and $A_3$ be 
three
disjoint annuli such that they contain each other in their outer exteriors, and 
for $i=1, 2, 3$, $A_i$ contains $x_i$ in 
its inner exterior. 
Let $B$ be an annulus that contains both $A_1$ and $A_2$ in its 
inner exterior, and $A_3$ in its outer exterior. 
Finally, let $C$ be an annulus that contains both $A_3$ and $B$ in its 
inner exterior. 
Now we define six one-way infinite curves, $\gamma_i^j$, 
for  $i=1, 2, 3$, $j=1, 2$.
For any fixed $i$, $i=1, 2, 3$, let $\gamma_i^1$ and $\gamma_i^2$ start
very close to $x_i$ and form a $k/4$-spiral in $A_i$.
In the outer exterior of $A_i$, let $\gamma_i^1$ and $\gamma_i^2$ form a
cable $\gamma_i$. Let $\gamma_1$ and $\gamma_2$ form a $k/4$-spiral in $B$.
In the outer exterior of $B$, let $\gamma_1$ and   $\gamma_2$ form a
cable $\gamma$. Finally, let $\gamma$ and $\gamma_3$ form a $k/4$-spiral in $C$.
In the outer exterior of $C$ all six curves go to infinity. 

Now replace each  one-way infinite curve $\gamma_i^j$ by two 
one-way infinite curves with the same endpoint, so that they go very close to
each other. Each of these pairs of curves form a bi-infinite curve $\Gamma_i^j$, and
any two intersect at most $k$ times. For the rest of the proof we call them
{\em pseudolines}.

Let $\rho$ be a bi-infinite curve containing $x_1$, $x_2$ and $x_3$. 
Then $\rho$ contains 
at least two
transversals of each of $A_1$, $A_2$, $A_3$, $B$ and $C$. 
This means, by Observation~\ref{obs_spiral},
that in each of $A_1$, $A_2$ and $A_3$,
there is a pseudoline, say, $\Gamma_1^1$, $\Gamma_2^1$ and $\Gamma_3^1$, respectively, that 
crosses $\rho$ at least $k/2-2$ times. Moreover, $\rho$ crosses
one of the two cables in $B$ at least $k/4-1$ times, which implies that
$\rho$ crosses one of the pseudolines
$\Gamma_1^1$ or $\Gamma_2^1$, say, $\Gamma_1^1$, at least $k/2-2$
times in $B$.
Finally, $\rho$ crosses the two cables in $C$ together at least $k/2-2$ times. Hence, in $C$, the curve $\rho$ has at least 
$k/8 - 1/2$ crossings with $\gamma$, or at least $3(k/8 - 1/2)$ crossings with $\gamma_3$.
In the first case $\rho$ crosses 
$\Gamma_1^1$ at least $5k/4-5$ times, in the second case it crosses 
$\Gamma_3^1$ at least $5k/4-5$ times.

(ii) For the second part of the theorem, we iterate the construction from the
proof of part (i) $m$ times. 

Let $P(k,0)$ be the following arrangement. Take a point $x$ in the plane and 
an annulus $A$ around it (that is, $x$ is in the inner exterior of $A$).
Let $\gamma^1$ and $\gamma^2$ be two one-way infinite curves, both starting near $x$
and forming a $k/4$-spiral in $A$.

Suppose that we have already defined an arrangement $P(k,i)$ containing
$2^i$ points 
and $2^{i+1}$  one-way infinite curves.
Take two disjoint copies of $P(k,i)$, and an annulus $B$ that contains all annuli of 
both copies in its internal exterior. Merge 
all curves of each copy of $P(k,i)$ into a cable and let the two cables form a
$k/4$-spiral in $B$. The resulting arrangement is $P(k,i+1)$.

Once the arrangement $P(k,m)$ is constructed, take two copies of each curve in $P(k,m)$ and join their endpoints to form a bi-infinite curve, thus obtaining a $k$-pseudoline arrangement $P'(k,m)$. 
Let $X_m$ be the set of $2^m$ points in the centers of the innermost annuli of $P'(k,m)$. By induction, every bi-infinite curve
containing all the points of $X_m$ crosses some pseudoline of $P'(k,m)$ at least
$(m+1)(k/2-2)$ times.
\end{proof}


\section{Concluding Remarks}

Our lower bound in Theorem~\ref{linearbound} for $k>1$ is weaker than for $k=1$. The reason
is that for
$k>1$, we
could not prove that a saturated $k$-simple topological graph cannot contain
an isolated triangle. The main difficulty is that for $k>1$, a triangle can
cross itself, and our proof for Lemma~\ref{triangle} does not work in this
case.

\begin{problem} 
\begin{enumerate}
\item[(i)] Is there a saturated $k$-simple topological graph,
for some $k\ge 2$, that contains an isolated triangle?
\item[(ii)] Is there a {\em disconnected\/} saturated $k$-simple topological graph,
for some $k$?
\end{enumerate}
\end{problem}

Problem 1 (ii) is open for every $k\ge 1$.

\bigskip

It follows from Theorem~\ref{klconstruction} (ii) that there is no
$(k,l)$-saturated graph with $l>2k$. By Theorem~\ref{3k/2saturated},
there is a $(k,l)$-saturated graph if $l\le \lceil 3k/2\rceil$.

\begin{problem}
Is there a $(k,l)$-saturated graph with $k\ge 2$ and $l>\lceil 3k/2\rceil$?
\end{problem}

In Theorem~\ref{pseudoline} we have shown that for sufficiently large $k$, there is a $(3, k+1)$-forcing arrangement of $k$-pseudolines. On the other hand, it is easy to see that there are no $(1, k+1)$-forcing arrangements of $k$-pseudolines. 

\begin{problem}
Is there a $(2, k+1)$-forcing arrangement of $k$-pseudolines for some $k\ge 3$?
\end{problem}

We assumed that in a  $k$-simple topological graph, no edge can cross itself.
For any $k$, a graph drawn in the plane is called a 
{\em $k$-complicated topological graph\/} if any two edges have at most $k$ points
in common, and an edge is allowed to cross itself, at most $k$ times.
Somewhat surprisingly, for saturated $k$-complicated topological graphs we cannot even
prove Lemma~\ref{isol+deg1} part (ii). We can only prove that a saturated
$k$-complicated topological graph
does not have isolated vertices. Therefore, the best lower bound we have for
the 
minimum number of edges of a saturated
$k$-complicated topological graph is $c_k(n)\ge n/2$.
On the other hand, for $k\ge 6$, using self-crossings, we can improve our upper bound
constructions from the proof of Theorem~\ref{linearbound} to obtain that 
$c_k(n)\le 5n/2$. We sketch the construction here.

Suppose that $n$ is even, $k\ge 6$, and let
$m=n/2$.  The construction consists of $2m+1$ consecutive blocks, 
$A_0,B_1,\allowbreak A_1,\allowbreak B_2,\allowbreak
\ldots,\allowbreak B_{m},\allowbreak A_{m}$, 
in this order, from bottom to top. 

For every $i$, $1\le i\le m$, 
let $u_i$ be a vertex on the lower boundary of $A_{i-1}$ and let $v_i$ be a vertex on the lower boundary of $B_i$. Let $\alpha_i$ be an edge joining $u_i$ and
$v_i$. The block $B_i$ is a $k$-spiral, and both of its cables are formed by
$\alpha_i$.
For $2\le i\le m$, the block $A_i$ is a 
$3$-spiral,
one cable is
formed by $\alpha_i$, the other one is formed by a folded curve
$\alpha_{i-1}$ (that is, two intervals of $\alpha_{i-1}$). 
Any curve $\kappa$ from $A_{i-2}$ to $A_{i}$
has to cross $\alpha_{i-1}$ at least $k-1$ times in $B_{i-1}$, and $\alpha_i$ at least $k-1$ times in $B_i$. In $A_{i-1}$, the curve $\kappa$ also crosses one of the curves $\alpha_i$ or $\alpha_{i-1}$ 
at least 
twice, since $\alpha_{i-1}$ is folded in $A_{i-1}$. It follows that when we extend this graph to a saturated 
$k$-complicated topological graph, each vertex has degree at most five. 

Note that using $\lfloor k/2\rfloor$-spirals in place of the $3$-spirals, we
obtain a $(k,l)$-saturated $k$-complicated topological graph with $l=5k/3 -
O(1)$.




\section*{Acknowledgements}
The first author thanks Radoslav Fulek and J\'anos Bar\'at for bringing the topic to his attention.

\end{document}